\documentclass{article}

\usepackage{amsmath}
\usepackage{amssymb}
\usepackage{amsthm}
\usepackage{amscd}

\newtheorem{thm}{Theorem}[section]
\newtheorem{prop}[thm]{Proposition}
\newtheorem{lem}[thm]{Lemma}
\newtheorem{cor}[thm]{Corollary}

\theoremstyle{definition}
\newtheorem{dfn}[thm]{Definition}

\theoremstyle{remark}
\newtheorem{rem}{Remark}
\newtheorem*{acknowledgment}{Acknowledgments}

\newcommand{\C}{\mathbb{C}}
\newcommand{\N}{\mathbb{N}}

\newcommand{\R}{\mathbb{R}}
\newcommand{\T}{\mathbb{T}}
\newcommand{\Z}{\mathbb{Z}}

\newcommand{\bH}{\mathbb{H}}

\newcommand{\G}{\mathcal{G}}
\renewcommand{\H}{\mathcal{H}}
\newcommand{\E}{\mathcal{E}}

\newcommand{\id}{\mathrm{id}}
\newcommand{\End}{\mathrm{End}}
\newcommand{\Ker}{\mathrm{Ker}}
\newcommand{\Hom}{\mathrm{Hom}}

\newcommand{\Cech}{$\check{\textrm{C}}$ech {}}
\newcommand{\Poincare}{Poincar\'e {}}
\newcommand{\im}{i}
\renewcommand{\l}{\ell}
\renewcommand{\u}{\underline}

\def\til#1{ \tilde{#1} }
\def\SDC#1{ \Z({#1})^\infty_D }

\def\abs#1{ \lvert{#1}\rvert }
\def\p#1{ \lVert{#1}\rVert }

\title{Projective unitary representations of
smooth Deligne cohomology groups}

\author{Kiyonori Gomi
\thanks{The author's research is supported by 
Research Fellowship of the Japan Society 
for the Promotion of Science for Young Scientists.}}

\date{}

\begin{document}

\maketitle

\begin{abstract}
We construct projective unitary representations of the smooth Deligne cohomology group of a compact oriented Riemannian manifold of dimension $4k+1$, generalizing positive energy representations of the loop group of the circle. We also classify such representations under a certain condition. The number of the equivalence classes of irreducible representations is finite, and is determined by the cohomology of the manifold.
\end{abstract}

\tableofcontents


\section{Introduction}
\label{sec:introduction}

The \textit{loop group} of a compact Lie group and its projective representations, called the \textit{positive energy representations} \cite{P-S}, have been studied involving various areas of mathematics. They are of physicist's interest as well, due to applications to quantum field theory.

From the viewpoint of ``higher abelian gerbes,'' a central extension of a certain smooth Deligne cohomology group is introduced in \cite{Go}. This is a generalization of a central extension of the loop group $L\T = C^\infty(S^1, \T)$ of the unit circle $\T = \{ z \in \Z |\ \abs{z} = 1 \}$. By definition, the smooth Deligne cohomology $H^q(M, \SDC{p})$ of a smooth manifold $M$ is the hypercohomology of the complex of sheaves:
$$
\SDC{p} : 
\Z \longrightarrow
\u{A}^0_M \overset{d}{\longrightarrow}
\u{A}^1_M \overset{d}{\longrightarrow}
\cdots \overset{d}{\longrightarrow}
\u{A}^{p-1}_M \longrightarrow
0 \longrightarrow \cdots,
$$
where $\Z$ denotes the constant sheaf, and $\u{A}^n_M$ the sheaf of germs of $\R$-valued differential $n$-forms on $M$. Suppose that $M$ is compact, oriented and $(4k+1)$-dimensional. For brevity, we put $\G(M) = H^{2k+1}(M, \SDC{2k+1})$. Then the cup product and the integration for smooth Deligne cohomology yield a non-trivial natural group 2-cocycle $S_M : \G(M) \times \G(M) \to \R/\Z$. Hence we have the corresponding central extension $\til{\G}(M)$ of $\G(M)$:
$$
\begin{CD}
1 @>>>
\T @>>>
\til{\G}(M) @>>>
\G(M) @>>>
1.
\end{CD}
$$
For instance, we consider the case of $k = 0$ and $M = S^1$. Then $\G(S^1)$ is isomorphic to the loop group $L\T$, and $\til{\G}(S^1)$ the central extension $\widehat{L\T}/\Z_2$, where $\widehat{L\T}$ is the \textit{universal central extension} \cite{P-S} of $L\T$. 

\bigskip

The purpose of this paper is to classify representations of $\til{\G}(M)$, that is, projective representations of $\G(M)$ with their cocycle $e^{2\pi\im S_M}$, under a certain condition. To state the main result precisely, we introduce some notions and notations. First of all, we remark that the group $\G(M) = H^{2k+1}(M, \SDC{2k+1})$ fits into the exact sequence:
\begin{equation} \label{iformula:exact_sequence}
0 \to 
A^{2k}(M)/A^{2k}(M)_\Z \to 
H^{2k+1}(M, \SDC{2k+1}) \to 
H^{2k+1}(M, \Z) \to 
0,
\end{equation}
where $A^{2k}(M)$ denotes the group of $2k$-forms on $M$, and $A^{2k}(M)_\Z$ the subgroup consisting of closed integral ones. We put $\G^0(M) = A^{2k}(M)/A^{2k}(M)_\Z$.

When $M$ is endowed with a Riemannian metric, we let $\bH^{2k}(M)$ be the group of harmonic $2k$-forms on $M$, and $\bH^{2k}(M) = \bH^{2k}(M) \cap A^{2k}(M)_\Z$ the subgroup consisting of integral ones. We denote by $\mathcal{X}(M)$ the group of homomorphisms $\lambda : \bH^{2k}(M)/\bH^{2k}(M)_\Z \to \R/\Z$.

\begin{lem} \label{ilem:identity_component}
Let $M$ be a compact oriented smooth $(4k+1)$-dimensional Riemannian manifold. For $\lambda \in \mathcal{X}(M)$, there exists an irreducible projective unitary representation $(\rho_\lambda, H_\lambda)$ of $\G^0(M)$ on a Hilbert space $H_\lambda$ with its cocycle $e^{2\pi\im S_M}$.
\end{lem}

In accordance with \cite{P-S}, the meaning of ``irreducible'' is that $H_\lambda$ contains no non-trivial invariant closed subspace. Notice that $(\rho_\lambda, H_\lambda)$, ($\lambda \in \mathcal{X}(M)$) are inequivalent representations. (We mean by an ``equivalence'' a continuous linear isomorphism compatible with the actions.)

\begin{dfn} \label{idfn:admissible}
An \textit{admissible representation} of $\G(M)$ is defined to be a projective unitary representation $(\rho, \H)$ of $\G(M)$ on a Hilbert space $\H$ with its cocycle $e^{2\pi\im S_M}$ such that: 
\begin{list}{}{\parsep=-2pt\topsep=4pt}
\item[(i)]
there is a map $m : \mathcal{X}(M) \to \Z_{\ge 0}$; and

\item[(ii)]
there is an equivalence of projective representations of $\G^0(M)$:
$$
\theta : \
\widehat{\bigoplus}_{\lambda \in \mathcal{X}(M)} 
\bigoplus^{m(\lambda)} H_\lambda
\longrightarrow
\H|_{\G^0(M)}.
$$
\end{list}
\end{dfn}

In the definition above, $\widehat{\oplus}$ stands for the Hilbert space direct sum. (We notice that Definition \ref{idfn:admissible} is motivated by \cite{F-H-T}. It may be possible to reformulate the admissibility above in terms of representations of the Lie algebra of $\G(M)$.)

Now, our main result is:

\begin{thm} \label{ithm:classification}
Let $M$ be a compact oriented $(4k+1)$-dimensional Riemannian manifold. 

(a) Any admissible representation of $\G(M)$ is equivalent to a finite direct sum of irreducible admissible representations.

(b) The number of the equivalence classes of irreducible admissible representations of $\G(M)$ is $2^br$, where $b = b_{2k+1}(M)$ is the Betti number, and $r$ is the number of the elements of the finite set $\{ t \in H^{2k+1}(M, \Z) |\ 2t = 0 \}$.
\end{thm}

For example, in the case of $k = 0$ and $M = S^1$, we have $\G(S^1) \cong L\T$ as mentioned. Since $H^1(S^1, \Z) \cong \Z$, the number of the equivalence classes of irreducible admissible representations of $\G(S^1)$ is 2, which coincides with that of irreducible positive energy representations of $L\T$ of level 2, (\cite{P-S,S}). We notice that any irreducible admissible representation of $\G(S^1)$ gives rise to an irreducible positive energy representation of $L\T$ of level 2, and vice verse, because our construction of the former coincides with that of the latter described in \cite{P-S}.

As another example, we consider the case of $k = 1$. Let $M$ be a compact oriented 5-dimensional Riemannian manifold whose second and third integral cohomology groups are torsion free. Then the number of the equivalence classes of irreducible admissible representations of $\G(M)$ is $2^b$. On the other hand, Henningson \cite{He} obtained $2^b$ Hilbert spaces in the quantum theory of a chiral 2-form on the 6-dimensional spacetime $\R \times M$. We can find a natural one to one correspondence between the parameterization of the quantum Hilbert spaces and that of the irreducible representations of $\G(M)$. In addition, using the basis of each quantum Hilbert space given in \cite{He}, we can identify it with the underlying Hilbert space of the corresponding representation in a natural way.

\bigskip

As generalizations of some properties of positive energy representations of the loop group $L\T$, the following result is also proved in this paper:

\begin{thm} \label{ithm:property}
Let $M$ be a compact oriented $(4k+1)$-dimensional Riemannian manifold, and $(\rho, \H)$ an admissible representation of $\G(M)$.

(a) $(\rho, \H)$ is continuous.

(b) There is an invariant dense subspace $\E \subset \H$ such that:
\begin{list}{}{\parsep=-2pt\topsep=4pt}
\item[(i)]
$(\rho, \E)$ extends to a projective representation of $\G(M)_\C$;

\item[(ii)]
the map $\rho(\cdot)\phi : \G(M)_\C \to \E$ is continuous for $\phi \in \E$.
\end{list}
\end{thm}

Some remarks are in order. In Theorem \ref{ithm:property} (a), the group $\G(M)$ is regarded as a topological group. The topology on $\G(M)$ is induced by the Riemannian metric on $M$ so that the group 2-cocycle $S_M : \G(M) \times \G(M) \to \R/\Z$ gives rise to a continuous map. (See Section \ref{sec:topology} for detail.) We mean by a ``continuous representation'' a representation $(\rho, \H)$ such that $\rho : \G(M) \times \H \to \H$ is a continuous map. In Theorem \ref{ithm:property} (b), the group $\G(M)_\C$ is the $(2k+1)$th hypercohomology of the complex of sheaves
$$
\Z \longrightarrow
\u{A}^0_{M, \C} \overset{d}{\longrightarrow}
\u{A}^1_{M, \C} \overset{d}{\longrightarrow}
\cdots \overset{d}{\longrightarrow}
\u{A}^{2k}_{M, \C} \longrightarrow
0 \longrightarrow \cdots,
$$
where $\u{A}^q_{M, \C}$ is the sheaf of germs of $\C$-valued $q$-forms on $M$. We can think of $\G(M)_\C$ as a ``complexification'' of $\G(M)$. In fact, we have $\G(S^1)_\C \cong L\C^*$, where $\C^* = \C - \{ 0 \}$. Notice that, for $f \in \G(M)_\C$, the linear map $\rho(f) : \E \to \E$ is unbounded in general.

\medskip 

A motivation of this paper comes from Wess-Zumino-Witten models. In this field, a remarkable result of Freed, Hopkins and Teleman \cite{F-H-T} relates an equivariant twisted $K$-theory to the {Verlinde algebra}, which is as a module generated by the equivalence classes of irreducible positive energy representations of a loop group. It seems interesting to pursue an analogous relation between the projective representations of $\G(M)$ and an equivariant twisted $K$-theory.

\bigskip

The organization of the present paper is as follows. 

Section \ref{sec:Deligne_coh} -- \ref{sec:harmonic_splitting} are preliminaries. In Section \ref{sec:Deligne_coh}, we recall some basic facts on smooth Deligne cohomology. We also introduce the group 2-cocycle $S_M$ and the central extension $\til{\G}(M)$. In Section \ref{sec:topology}, we introduce an inner product on the space of differential forms $A^{2k}(M)$. We utilize it to topologize $\G(M)$. The inner product is also a key to constructing projective representations of $\G(M)$. In Section \ref{sec:harmonic_splitting}, we study in detail a certain splitting of the exact sequence (\ref{iformula:exact_sequence}), which leads to a decomposition of $\G(M)$.

Then, in Section \ref{sec:construction}, we construct projective representations of $\G(M)$, generalizing that of positive energy representations of $L\T$ described in \cite{P-S}. Section \ref{sec:classification} is the heart of this paper, in which Theorem \ref{ithm:classification} is proved. Section \ref{sec:properties} is devoted to the proof of Theorem \ref{ithm:property}. 

Finally, in Appendix, we prove some general properties of unitary representations on Hilbert spaces. A kind of ``Schur's lemma'' shown here plays the most important role in the classification. In their proof, we use some tools in functional analysis, such as the spectral decomposition theorem. To make the point of an argument transparent, we separated it from the main text.

\bigskip

From Section \ref{sec:topology} through Section \ref{sec:properties}, we let $k$ denote a non-negative integer, and $M$ a compact oriented smooth $(4k+1)$-dimensional Riemannian manifold without boundary.


\section{Smooth Deligne cohomology}
\label{sec:Deligne_coh}

We review some basic facts on smooth Deligne cohomology \cite{Bry,D-F,E-V}. Then we introduce the group cocycle $S_M$ by means of the cup product and the integration for smooth Deligne cohomology.


\subsection{Definition of smooth Deligne cohomology}

Let $M$ be a (finite dimensional) smooth manifold. For a non-negative integer $q$, we denote by $\u{A}^q_M$ the sheaf of germs of $\R$-valued differential $q$-forms on $M$.

\begin{dfn}
For a non-negative integer $p$, we define the \textit{smooth Deligne complex} $\SDC{p}$ to be the following complex of sheaves:
$$
\SDC{p} : \
\Z \longrightarrow
\u{A}^0_M \overset{d}{\longrightarrow}
\u{A}^1_M \overset{d}{\longrightarrow}
\cdots \overset{d}{\longrightarrow}
\u{A}^{p-1}_M \longrightarrow
0 \longrightarrow \cdots,
$$
where $\Z$ is the constant sheaf and is located at degree $0$ in the complex. The \textit{smooth Deligne cohomology group} $H^q(M, \SDC{p})$ of $M$ is defined to be the hypercohomology of $\SDC{p}$.
\end{dfn}

We can readily see the natural isomorphism $H^1(M, \SDC{1}) \cong C^\infty(M, \T)$, where we denote by $\T = \{ z \in \C |\ \abs{z} = 1 \}$ the unit circle. 

\medskip

The following proposition, which is easily shown, reveals the relation between the smooth Deligne cohomology and other cohomology groups.

\begin{prop}[\cite{Bry}] \label{prop:Deligne_coh_exact_seq}
Let $p$ be a positive integer.

(a) If $0 \le q < p$, then $H^q(M, \SDC{p}) \cong H^{q-1}(M, \R/\Z)$.

(b) If $p = q$, then $H^p(M, \SDC{p})$ fits into the following exact sequences:
\begin{gather}
0 \to
H^{p-1}(M, \R/\Z) \to
H^p(M, \SDC{p}) \overset{\delta}{\to}
A^p(M)_{\Z} \to 0, \\
0 \to
A^{p-1}(M) / A^{p-1}(M)_{\Z} \overset{\iota}{\to}
H^p(M, \SDC{p}) \overset{\chi}{\to}
H^p(M, \Z) \to 0,
\label{formula:second_exact_seq}
\end{gather}
where $A^q(M)_{\Z} \subset A^q(M)$ is the subgroup of closed integral $q$-forms on $M$. 

(c) If $p < q$, then $H^q(M, \SDC{p}) \cong H^q(M, \Z)$.
\end{prop}

We have the following compositions of homomorphisms in Proposition \ref{prop:Deligne_coh_exact_seq}:
\begin{align*}
A^{p-1}(M)/A^{p-1}(M)_\Z \overset{\iota}{\to} 
H^p(M, \SDC{p}) \overset{\delta}{\to} 
A^p(M), \\
H^{p-1}(M, \R/\Z) \to 
H^p(M, \SDC{p}) \overset{\chi}{\to} 
H^p(M, \Z).
\end{align*}
The first composition coincides with the homomorphism induced by the exterior derivative $d : A^{p-1}(M) \to A^p(M)$. The second composition coincides with the connecting homomorphism $\beta : H^{p-1}(M, \R/\Z) \to H^p(M, \Z)$ in the exact sequence of cohomology groups induced by $\Z \to \R \to \R/\Z$.

\smallskip

Let $H^p_{DR}(M)$ be the de Rham cohomology group. The operation of taking the de Rham cohomology class of a closed form induces the natural homomorphism $A^p(M)_\Z \to H^p_{DR}(M)$. The inclusion $\Z \to \R$ also induces the natural homomorphism $H^p(M, \Z) \to H^p(M, \R)$. By de Rham's theorem, we have the natural isomorphism $H^p(M, \R) \cong H^p_{DR}(M)$. These natural homomorphisms fit into the following commutative diagram:
\begin{equation}
\begin{CD}
H^p(M, \SDC{p}) @>{\delta}>> A^p(M)_\Z  @>>> H^p_{DR}(M) \\
@|                           @.              @VV{\cong}V \\
H^p(M, \SDC{p}) @>>{\chi}>   H^p(M, \Z) @>>> H^p(M, \R).
\label{formula:de_Rham_real_coefficients}
\end{CD}
\end{equation}


\subsection{Cup product and integration}

We recall the cup product and the integration for smooth Deligne cohomology.

\medskip

The \textit{cup product} is the natural homomorphism
$$
\cup : H^p(M, \SDC{p}) \otimes_\Z H^q(M, \SDC{q}) \to H^{p+q}(M, \SDC{p+q})
$$
induced from a homomorphism of sheaves of complexes. (We refer the reader to \cite{Bry,E-V} for the detail of the definition.) It is known that the cup product is associative and (graded) commutative:
$$
(f_p \cup f_q) \cup f_q = f_p \cup (f_q \cup f_r), \quad
f_p \cup f_q = (-1)^{pq} f_q \cup f_p,
$$
where $f_i \in H^i(M, \SDC{i})$ for $i = p, q, r$. Through the homomorphisms $\delta$ and $\chi$, the cup product on $H^p(M, \SDC{p})$ is compatible with the wedge product on $A^p(M)_\Z$ and with the cup product on $H^p(M, \Z)$:
$$
\delta(f \cup g) = \delta(f) \wedge \delta(g), \quad
\chi(f \cup g) = \chi(f) \cup \chi(g).
$$
If we restrict the cup product on $H^p(M, \SDC{p})$ to the subgroup $H^{p-1}(M, \R/\Z)$, then we obtain the homomorphism
\begin{equation}
H^{p-1}(M, \R/\Z) \otimes_\Z H^{q-1}(M, \R/\Z) \to
H^{p+q-1}(M, \R/\Z). \label{formula:restrict_cup_product_1}
\end{equation}
This is identified with the composition of the following homomorphisms:
\begin{align*}
H^{p-1}(M, \R/\Z) \otimes_\Z H^{q-1}(M, \R/\Z)
& \overset{\beta \otimes \id}{\longrightarrow}
H^{p}(M, \Z) \otimes_\Z H^{q-1}(M, \R/\Z), \\
H^{p}(M, \Z) \otimes_\Z H^{q-1}(M, \R/\Z)
& \overset{\cup}{\longrightarrow}
H^{p + q -1}(M, \R/\Z),
\end{align*}
where $\cup$ is the cup product induced by $\Z \otimes_\Z \R/\Z \to \R/\Z$. Similarly, if we restrict the cup product on $H^p(M, \SDC{p})$ to $A^{p-1}(M)/A^{p-1}(M)_\Z$, then we have the homomorphism
$$
(A^{p-1}(M)/A^{p-1}(M)_\Z) \otimes_\Z (A^{q-1}(M)/A^{q-1}(M)_\Z) \to
A^{p+q-1}(M)/A^{p+q-1}(M)_\Z.
$$
This homomorphism is given by $(\alpha, \beta) \mapsto \alpha \wedge d\beta$.

\bigskip

When $M$ is a compact oriented $m$-dimensional smooth manifold without boundary, the \textit{integration} for smooth Deligne cohomology is given as a natural homomorphism
$$
\int_M : \ H^{m+1}(M, \SDC{m+1}) \longrightarrow \R/\Z
$$
inducing the isomorphisms $H^{m+1}(M, \SDC{m+1}) \cong H^m(M, \R/\Z) \cong \R/\Z$ and $H^{m+1}(M, \SDC{m+1}) \cong A^m(M)/A^m(M)_\Z \cong \R/\Z$. (An explicit formula on the level of \Cech cocycles can be found in \cite{Go,G-T}.) The integration satisfies
\begin{equation}
\int_M \iota(\alpha) \cup f = \int_M \alpha \wedge \delta(f) \mod \Z,
\label{formula:integration_Deligne}
\end{equation}
where $\alpha \in A^p(M)/A^p(M)_\Z$, $f \in H^q(M, \SDC{q})$ and $p + q = m$.


\subsection{Central extension}

Let $n$ be a non-negative integer. For a compact oriented smooth $(2n+1)$-dimensional smooth manifold $M$ without boundary, we denote by $\G(M)$ the smooth Deligne cohomology group $\G(M) = H^{n+1}(M, \SDC{n+1})$. 

\begin{dfn}[\cite{Go}]
We define the group 2-cocycle $S_M : \G(M) \times \G(M) \to \R/\Z$ by $S_M(f, g) = \int_M f \cup g$. We also define the group $\til{\G}(M)$ by the set $\til{\G}(M) = \G(M) \times \T$ endowed with the group multiplication
$$
(f, u) \cdot (g, v) = (f + g, uv \exp2\pi\im S_M(f, g)).
$$
\end{dfn}

By definition, $\til{\G}(M)$ gives a central extension of $\G(M)$ by $\T$:
\begin{equation}
\begin{CD}
1  @>>>
\T @>>>
\til{\G}(M) @>>>
\G(M) @>>>
1.
\end{CD}
\label{formula:central_extension}
\end{equation}

\begin{prop}[\cite{Go}]
As a central extension, $\til{\G}(M)$ is non-trivial, if and only if $n$ is even.
\end{prop}

Because of this result, the case of $n$ even is of our interest.

As is noticed, there is a natural isomorphism $H^1(M, \SDC{1}) \cong C^\infty(M, \R/\Z)$. In the case where $n = 0$ and $M = S^1$, we can identify $\G(S^1)$ with the loop group of the unit circle: $\G(S^1) \cong L\T$. It is known \cite{Bry-M2} that, under this identification, we have an isomorphism of central extensions $\til{\G}(S^1) \cong \widehat{L\T}/\Z_2$, where $\widehat{L\T}$ is the \textit{universal central extension} of the loop group $L\T$, (\cite{P-S,S}).


\section{Topology on the central extension}
\label{sec:topology}

In this section, we topologize the smooth Deligne cohomology $\G(M)$ so that the group 2-cocycle gives rise to a continuous map. 

\medskip

From this section to the end, we let $k$ be a non-negative integer, and $M$ a compact oriented smooth $(4k+1)$-dimensional Riemannian manifold without boundary. We denote by $\G(M)$ the smooth Deligne cohomology group $\G(M) = H^{2k+1}(M, \SDC{2k+1})$.

\subsection{Topology on the space of differential forms}

To make the smooth Deligne cohomology groups into topological groups, we begin with the following proposition regarding the space of differential forms.

\begin{prop} \label{prop:topology_differential_forms}
There exists a positive definite inner product 
$$
( \ , \ )_V : A^{2k}(M) \times A^{2k}(M) \to \R
$$ 
satisfying the following properties:

(a) If we topologize $A^{2k}(M)$ by $( \ , \ )_V$, then $A^{2k}(M)$ decomposes into mutually orthogonal closed subspaces with respect to $( \ , \ )_V$:
$$
A^{2k}(M) \cong d(A^{2k-1}(M)) \oplus \bH^{2k}(M) \oplus d^*(A^{2k+1}(M)),
$$
where $d^* : A^{2k+1}(M) \to A^{2k}(M)$ is defined to be $d^* = - * d *$ by means of the Hodge star operators, and $\bH^{2k}(M)$ is the space of Harmonic $2k$-forms on $M$.

(b) The Hilbert space $V$ given by the completion of $d^*(A^{2k+1}(M))$ with respect to $( \ , \ )_V$ is separable.

(c) There exists a linear map $J : V \to V$ such that:
\begin{list}{}{\parsep=-2pt\topsep=4pt}
\item[(i)] $J^2 = -1$;

\item[(ii)] $(Jv, Jv')_V = (v, v')_V$ for $v, v' \in V$;

\item[(iii)] $(\alpha, J\beta)_V = \int_M \alpha \wedge d\beta$ for $\alpha, \beta \in d^*(A^{2k+1}(M)) \subset V$.
\end{list}
\end{prop}

The crucial property is Proposition \ref{prop:topology_differential_forms} (c) (iii). One can realize $( \ , \ )_V$ as the inner product induced by the Sobolev $H^s$-norm with $s = 1/2$, and $J : V \to V$ as $J = \til{J}/\abs{\til{J}}$, where $\til{J} : d^*(A^{2k+1}(M)) \to d^*(A^{2k+1}(M))$ is the differential operator $\til{J} = * d$. In the proof below, we give a concrete realization using some basic facts about eigenforms of the Laplacian. (See \cite{Wa}, for example.)

\begin{proof}
We only describe constructions of $( \ , \ )_V$ and $J$, for we can verify their properties in an elementary way. The Riemannian metric on $M$ defines the Laplacian $\Delta = dd^* + d^*d$ and the $L^2$-inner product $( \alpha , \beta )_{L^2} = \int_M \alpha \wedge * \beta$ on $A^{2k}(M)$. Let $0 \le \l_1 \le \l_2 \le \cdots$ be the eigenvalues of $\Delta$, where each eigenvalue is included as many times as the dimension of its eigenspace. Let $\{ \psi_i \}_{i \in \N}$ be a sequence of corresponding orthonormalized eigenforms: $\Delta \psi_i = \l_i \psi_i$ and $(\psi_i, \psi_j)_{L^2} = \delta_{ij}$. Because of the Hodge decomposition theorem, we can take the eigenforms so that each $\psi_i$ belongs to either $d(A^{2k-1}(M))$, $\bH^{2k}(M)$ or $d^*(A^{2k+1}(M))$. We define the norm $\p{ \ \ }_V : A^{2k+1}(M) \to \R$ by
$$
\p{\alpha}_V^2 = \sum_{i = 1}^b(\alpha, \psi_i)_{L^2}^2 + \sum_{i = b + 1}^\infty \sqrt{\l_i}(\alpha, \psi_i)_{L^2}^2,
$$
where $b = b_{2k}(M)$. We notice $\sqrt{\l_i} < 1 + \l_i^2$, so that $\p{\alpha}_V^2 \le \p{\alpha}_{L^2}^2 + \p{\Delta \alpha}_{L^2}^2$. The norm $\p{ \ \ }_V$ induces the inner product $( \ , \ )_V$. 

To define $J$, we let $\{ \psi_{i(j)} \}_{j \in \N}$ be the subsequence consisting of the eigenforms belonging to $d^*(A^{2k+1}(M))$, and $\{ \l_{i(j)} \}_{j \in \N}$ the corresponding sequence of eigenvalues. If we put $\lambda_j = \l_{i(j)}$ and $\varphi_j = \psi_{i(j)}/\sqrt[4]{\lambda_j}$, then $\{ \varphi_j \}_{j \in \N}$ forms a complete orthonormal basis of $V$. Notice that $\tilde{J} = * d = - d^* *$ obeys $\tilde{J}^2 = - \Delta$ on $d^*(A^{2k+1}(M))$. We define the isometric map $J : d^*(A^{2k+1}(M)) \to V$ by
$$
J \alpha = 
\sum_{j = 1}^\infty 
(\alpha, \varphi_j)_V \frac{1}{\sqrt{\lambda_j}} \tilde{J} \varphi_j.
$$
The extension is $J : V \to V$. We remark that $( \ , \ )_V$ and $J$ are independent of the choice of eigenforms of $\Delta$, since each eigenspace of $\Delta$ is finite dimensional.
\end{proof}

When we make $A^{2k}(M)$ into a Hausdorff locally convex topological vector space by means of the inner product $( \ , \ )_V$, Proposition \ref{prop:topology_differential_forms} (a) leads to:

\begin{cor}
We have the following isomorphism of topological groups:
\begin{align*}
A^{2k}(M)_\Z 
&\cong 
d(A^{2k-1}(M)) \times \bH^{2k}(M)_\Z, \\
A^{2k}(M)/A^{2k}(M)_\Z 
&\cong 
(\bH^{2k}(M)/\bH^{2k}(M)_\Z) \times d^*(A^{2k+1}(M)),
\end{align*}
where $\bH^{2k}(M)_\Z = A^{2k}(M)_\Z \cap \bH^{2k}(M)$ is endowed with the induced topology.
\end{cor}

As is well-known, $\bH^{2k}(M)$ is a finite dimensional vector space. So the induced topology on $\bH^{2k}(M)$ coincides with the standard topology on $\R^b$, and $\bH^{2k}(M)_\Z$ is isomorphic to $\Z^b$ with the discrete topology, where $b = b_{2k}(M)$ is the $2k$th Betti number of $M$.

\subsection{Topology on the smooth Deligne cohomology}

Recall the exact sequence (\ref{formula:second_exact_seq}) in Proposition \ref{prop:Deligne_coh_exact_seq}:
\begin{equation}
0 \to
A^{2k}(M) / A^{2k}(M)_{\Z} \overset{\iota}{\to}
\G(M) \overset{\chi}{\to}
H^{2k+1}(M, \Z) \to 0.
\label{formula:second_exact_seq_again}
\end{equation}

\begin{lem}
There exists a splitting $\sigma : H^{2k+1}(M, \Z) \to \G(M)$ of the exact sequence (\ref{formula:second_exact_seq_again}). Hence we have the isomorphism of abelian groups
\begin{equation}
I_\sigma : \ 
(A^{2k}(M)/A^{2k}(M)_\Z) \times H^{2k+1}(M, \Z) \longrightarrow \G(M)
\label{fomula:splitting}
\end{equation}
given by $I_\sigma(\alpha, c) = \alpha + \sigma(c)$.
\end{lem}

\begin{proof}
Because $M$ is compact, $H^{2k+1}(M, \Z)$ is finitely generated. It is easy to see that $A^{2k}(M)/A^{2k}(M)_\Z$ is divisible, so that the exact sequence (\ref{formula:second_exact_seq_again}) splits.
\end{proof}

We make $(A^{2k}(M)/A^{2k}(M)_\Z) \times H^{n+1}(M, \Z)$ into a Hausdorff topological group by the topology on $A^{2k}(M)/A^{2k}(M)_\Z$ given in the previous subsection and by the discrete topology on $H^{2k+1}(M, \Z)$.

\begin{dfn} \label{dfn:topology_Deligne_coh}
Taking a splitting $\sigma$ of (\ref{formula:second_exact_seq_again}), we make $\G(M)$ into a Hausdorff topological group so that $I_\sigma$ is a homeomorphism. 
\end{dfn}

\begin{lem}
The topology on $\G(M)$ is independent of the choice of $\sigma$.
\end{lem}

\begin{proof}
Let $\sigma'$ be the other splitting. It suffices to prove that 
$$
I_{\sigma}^{-1} \circ I_{\sigma'} : \
(A^{2k}/A^{2k}_\Z) \times H^{2k+1}(M, \Z) \to
(A^{2k}/A^{2k}_\Z) \times H^{2k+1}(M, \Z)
$$ 
is continuous. We can readily see $I_{\sigma}^{-1} \circ I_{\sigma'}(\alpha, c) = (\alpha + \tau(c), c)$, where the homomorphism $\tau : H^{2k+1}(M, \Z) \to A^{2k}(M)/A^{2k}(M)_\Z$ is defined by the difference $\tau = \sigma' - \sigma$. Since the topology on $H^{2k+1}(M, \Z)$ is discrete, the homomorphism $\tau$ is clearly continuous. Hence $I_{\sigma}^{-1} \circ I_{\sigma'}$ is continuous as well.
\end{proof}

It is clear by definition that the identity component of $\G(M)$ is the subgroup $A^{2k}(M)/A^{2k}(M)_\Z$. From now on, we put $\G^0(M) = A^{2k}(M)/A^{2k}(M)_\Z$.

\subsection{Topology on the central extension}

We begin with introducing a useful notion:

\begin{dfn} \label{dfn:harmonic_splitting}
We define a \textit{harmonic splitting} $\sigma : H^{2k+1}(M, \Z) \to \G(M)$ to be a splitting of (\ref{formula:second_exact_seq_again}) such that $\delta(\sigma(c)) \in A^{2k+1}(M)$ is a harmonic form for each $c \in H^{2k+1}(M, \Z)$.
\end{dfn}

\begin{lem}
(a) There exists a harmonic splitting $\sigma : H^{2k+1}(M, \Z) \to \G(M)$. 

(b) For harmonic splittings $\sigma$ and $\sigma'$, there uniquely exists a homomorphism $\tau : H^{2k+1}(M, \Z) \to \bH^{2k}(M)/\bH^{2k}(M)_\Z$ such that $\sigma' = \sigma + \tau$.
\end{lem}

\begin{proof}
For (a), we construct a harmonic splitting $\sigma'$ from a given splitting $\sigma$. First, we take a basis of $H^{2k+1}(M, \Z)$ so that we have
$$
H^{2k+1}(M, \Z) \cong 
\Z e_1 \oplus \cdots \oplus \Z e_b \oplus 
\Z_{p_1} t_1 \oplus \cdots \oplus \Z_{p_\l} t_\l,
$$
where $e_1, \ldots, e_b$ generate the free part, and $t_1, \ldots t_\l$ the torsion part. Note that $\delta(\sigma(t_i)) = 0$. Let $\eta_i \in \bH^{2k+1}(M)$ be the unique harmonic form whose de Rham cohomology class coincides with that of $\delta(\sigma(e_i))$. So there is a $2k$-form $\alpha_i$ such that $\delta(\sigma(e_i)) + d \alpha_i = \eta_i$. Now we define the splitting $\sigma' : H^{2k+1}(M, \Z) \to \G(M)$ to be the linear extension of $\sigma'(e_i) = \sigma(e_i) + \alpha_i$ and $\sigma'(t_j) = \sigma(t_j)$. Then $\sigma'$ is harmonic. 

For (b), let $\tau : H^{2k+1}(M, \Z) \to A^{2k}(M)/A^{2k}(M)_\Z$ be the homomorphism defined by $\tau = \sigma' - \sigma$. The commutativity of the diagram (\ref{formula:de_Rham_real_coefficients}) implies that, for $c \in H^{2k+1}(M, \Z)$, the harmonic forms $\delta(\sigma(c))$ and $\delta(\sigma(c'))$ are identical, so that $d \tau(c) = 0$. Recall here the isomorphism of abelian groups:
$$
A^{2k}(M)/A^{2k}(M)_\Z \cong 
\bH^{2k}(M)/\bH^{2k}(M)_\Z \times d^*(A^{2k+1}(M)).
$$
According to this isomorphism, we write $\tau = \tau' + \tau''$, where $\tau' : H^{2k+1}(M, \Z) \to \bH^{2k}(M)/\bH^{2k}(M)_\Z$ and $\tau'' : H^{2k+1}(M, \Z) \to d^*(A^{2k+1}(M))$. Then $d\tau(c) = 0$ implies $d\tau''(c) = 0$. Because the only closed form contained in $d^*(A^{2k+1}(M))$ is the trivial one, we see $\tau''(c) = 0$. Hence $\tau = \tau'$, and the proof is completed.
\end{proof}

A choice of a splitting $\sigma$ induces the isomorphism of topological groups
$$
I_\sigma : 
(\bH^{2k}(M) / \bH^{2k}(M)_\Z) \times d^*(A^{2k+1}(M)) \times H^{2k+1}(M, \Z)
\longrightarrow \G(M),
$$
where $I_\sigma(\eta, \nu, c) = \eta + \nu + \sigma(c)$. If $\sigma$ is harmonic, then (\ref{formula:integration_Deligne}) gives:
\begin{multline}
S_M(I_\sigma(\eta, \nu, c), I_\sigma(\eta', \nu', c')) \\
=
\int_M \nu \wedge d\nu'
+
\int_M \left\{ 
\eta \wedge \delta(\sigma(c')) - \eta' \wedge \delta(\sigma(c))
\right\} 
+
(\sigma^*S_M)(c, c'),
\label{formula:expression_S_harmonic}
\end{multline}
where $(\sigma^*S_M)(c, c') = S_M(\sigma(c), \sigma(c'))$. Thanks to the harmonicity of $\sigma$, terms such as $\int_M \nu \wedge \delta (\sigma(c'))$ are absent in the expression above.

\begin{prop} \label{prop:continuity_2_cocycle}
$S_M : \G(M) \times \G(M) \to \R/\Z$ is continuous.
\end{prop}

\begin{proof}
By the help of the expression (\ref{formula:expression_S_harmonic}), the continuity of $S_M$ amounts to that of the following maps:
\begin{align*}
H^{2k+1}(M, \Z) &\times H^{2k+1}(M, \Z) \to \R/\Z, &
(c, c') &\mapsto \int_M \sigma(c) \cup \sigma(c'), \\
\bH^{2k}(M) &\times H^{2k+1}(M, \Z) \to \R, &
(\eta, c) &\mapsto \int_M \eta \wedge \delta(\sigma(c)), \\
d^*(A^{2k+1}(M)) &\times d^*(A^{2k+1}(M)) \to \R, &
(\nu, \nu') &\mapsto \int_M \nu \wedge d\nu'.
\end{align*}
Since the topology on $H^{2k+1}(M, \Z)$ is discrete, we can readily see that the first and the second map are continuous. The continuity of the third map follows from Proposition \ref{prop:topology_differential_forms} (c).
\end{proof}

\begin{cor}
We can make $\til{\G}(M)$ into a Hausdorff topological group so that (\ref{formula:central_extension}) is an exact sequence of Hausdorff topological groups.
\end{cor}

\begin{proof}
It suffices to topologize $\til{\G}(M) = \G(M) \times \T$ by using the topology on $\G(M)$ and the standard topology on $\T$.
\end{proof}

\begin{rem}
We can use a finer topology on $A^{2k}(M)$ to make $\G(M)$ into a topological group. The choice of the topology in the present paper is essential only when we combine Proposition \ref{prop:topology_differential_forms} with Proposition \ref{prop:Heisenberg} in constructing projective representation of $\G(M)$.
\end{rem}


\section{Harmonic splitting}
\label{sec:harmonic_splitting}

This section is devoted to a detailed study of harmonic splittings introduced in Definition \ref{dfn:harmonic_splitting}. The results obtained here will be used in both the construction and the classification of projective representations of $\G(M)$.

\subsection{Linking form}

We first see that the bilinear map
$$
\sigma^* S_M : \ H^{2k+1}(M, \Z) \times H^{2k+1}(M, \Z) \to \R/\Z
$$
given by taking a splitting $\sigma$ has some part independent of the choice of $\sigma$. 

We denote by $T^{2k+1} = \Ker\{ H^{2k+1}(M, \Z) \to H^{2k+1}(M, \R) \}$ the torsion subgroup in $H^{2k+1}(M, \Z)$.

\begin{lem}
Let $\sigma : H^{2k+1}(M, \Z) \to \G(M)$ be a splitting. The maps
\begin{align}
T^{2k+1} \times T^{2k+1} &\longrightarrow \R/\Z, &
(t, t') &\mapsto S_M(\sigma(t), \sigma(t')), 
\label{formula:linking_pairing} \\
H^{2k+1}(M, \Z) &\longrightarrow \R/\Z, &
c &\mapsto S_M(\sigma(c), \sigma(c))
\label{formula:Steenrod_square}
\end{align}
are independent of the choice of $\sigma$.
\end{lem}

\begin{proof}
Let $\sigma'$ be another splitting. By setting $\tau = \sigma' - \sigma$, we define the homomorphism $\tau : H^{2k+1}(M, \Z) \to A^{2k}(M)/A^{2k}(M)_\Z$. Then we have 
\begin{multline*}
S_M(\sigma'(t), \sigma'(t')) - S_M(\sigma(t), \sigma(t')) \\
=
\int_M \tau(t) \wedge \delta (\sigma(t')) 
- \int_M \tau(t') \wedge \delta (\sigma(t))
+ \int_M \tau(t) \wedge d (\tau(t')) \mod \Z.
\end{multline*}
Because $A^{2k+1}(M)$ is torsion free, we have $\delta \circ \sigma(t) = d \circ \tau(t) = 0$ for $t \in T^{2k+1}$. Hence (\ref{formula:linking_pairing}) is independent of the choice of the splitting. Similarly, we have
$$
S_M(\sigma'(c), \sigma'(c)) - S_M(\sigma(c), \sigma(c))
= \int_M \tau(c) \wedge d(\tau(c)) \mod \Z
$$
for $c \in H^{2k+1}(M, \Z)$. Because $\tau(c) \wedge d(\tau(c)) = \frac{1}{2} d (\tau(c) \wedge \tau(c))$, Stokes' theorem implies that (\ref{formula:Steenrod_square}) is also independent of the choice of $\sigma$.
\end{proof}

\begin{dfn}
We define a map $L_M : T^{2k+1} \times T^{2k+1} \to \R/\Z$ by (\ref{formula:linking_pairing}).
\end{dfn}

\begin{lem}
The map $L_M$ coincides with the \textit{linking form} (\cite{S-T}).
\end{lem}

\begin{proof}
Let $t, t' \in T^{2k+1}$ be torsion elements. Then there is $\til{t} \in H^{2k}(M, \R/\Z)$ such that $\beta(\til{t}) = t$, where $\beta : H^{2k}(M, \R/\Z) \to H^{2k+1}(M, \Z)$ is the connecting homomorphism. By virtue of the description of the cup product given in (\ref{formula:restrict_cup_product_1}), we obtain $L_M(t, t') = \til{t} \cup t'$, which leads to the present lemma.
\end{proof}

Notice that $2 \sigma^* S_M(c, c) = 0$ for $c \in H^{2k+1}(M, \Z)$ by the graded commutativity of the cup product for smooth Deligne cohomology.

\begin{rem}
The map (\ref{formula:Steenrod_square}) coincides with the composition:
$$
H^{2k+1}(M, \Z) \overset{\mod 2}{\longrightarrow} 
H^{2k+1}(M, \Z_2) \overset{\mathrm{Sq}^{2k}}{\longrightarrow}
H^{4k+1}(M, \Z_2) \cong (\frac{1}{2}\Z)/\Z \to \R/\Z,
$$
where $\mathrm{Sq}^{2k}$ is the Steenrod squaring operation (\cite{Spa}). This fact follows from the result in \cite{Go2}, because there are natural isomorphisms between the smooth Deligne cohomology groups and the groups of differential characters \cite{Bry}.
\end{rem}

\subsection{Compatible harmonic splitting}

Since $H^{2k+1}(M, \Z)$ is a finitely generated abelian group, we can decompose it into a free part and the torsion part. In this paper, we use the following terminology for convenience.

\begin{dfn}
A \textit{decomposition} of $H^{2k+1}(M, \Z)$ is an isomorphism of abelian groups $\omega : H^{2k+1}(M, \Z) \cong F^{2k+1} \oplus T^{2k+1}$, where $F^{2k+1} \subset H^{2k+1}(M, \Z)$ is a free subgroup of rank $b = b_{2k+1}(M)$.
\end{dfn}

A decomposition $\omega$ induces the homomorphism $\bar{\omega} : H^{2k+1}(M, \Z) \to T^{2k+1}$ such that $\bar{\omega}(t) = t$ for $t \in T^{2k+1}$. Conversely, such a homomorphism $\bar{\omega}$ recovers the decomposition $\omega : H^{2k+1}(M, \Z) \cong F^{2k+1} \oplus T^{2k+1}$ by $F^{2k+1} = \Ker \ \bar{\omega}$ and $\omega(c) = (c - \bar{\omega}(c), \bar{\omega}(c))$. This correspondence will be freely used.

\begin{dfn}
Let $\omega : H^{2k+1}(M, \Z) \cong F^{2k+1} \oplus T^{2k+1}$ be a decomposition given. A splitting $\sigma : H^{2k+1}(M, \Z) \to \G(M)$ is said to be \textit{compatible with $\omega$} when $\sigma^*S_M(t, e) = 0$ for all $t \in T^{2k+1}$ and $e \in F^{2k+1}$.
\end{dfn}

\begin{prop} \label{prop:compatible_harmonic_splitting}
Let $\omega : H^{2k+1}(M, \Z) \cong F^{2k+1} \oplus T^{2k+1}$ be a decomposition.

(a) There is a harmonic splitting compatible with $\omega$.

(b) If harmonic splittings $\sigma$ and $\sigma'$ are compatible with $\omega$, then $\sigma (t) = \sigma'(t)$ for $t \in T^{2k+1}$. Thus, in particular, Harmonic splittings compatible with $\omega$ are in one to one correspondence with homomorphisms $F^{2k+1} \to \bH^{2k}(M)/\bH^{2k}(M)_\Z$.
\end{prop}

\begin{proof}
Let $\{ e_1, \ldots, e_b \}$ be a basis of $F^{2k+1}$. Then there is an isomorphism
$$
H^{2k+1}(M, \Z) \cong 
\Z e_1 \oplus \cdots \oplus \Z e_b
\oplus
\Z_{p_1} t_1 \oplus \cdots \oplus \Z_{p_\l} t_\l,
$$
where $t_1, \ldots, t_\l$ generate the torsion part $T^{2k+1}$. By the \Poincare duality and the universal coefficients theorem (see \cite{Spa}), there are cohomology classes $h_1, \ldots, h_b \in H^{2k}(M, \Z)$ generating a free subgroup of rank $b$ such that $\langle h_i \cup e_j, [M] \rangle = \delta_{ij}$. We denote by $\eta_j \in \bH^{2k}(M)$ the harmonic representative of the image of $h_j$ under the natural homomorphism $H^{2k}(M, \Z) \to H^{2k}(M, \R) \cong H^{2k}_{DR}(M)$. 

Now, take and fix a harmonic splitting $\sigma_0$. For $i = 1, \ldots, \l$ and $j = 1, \ldots, b$, we take $Z_{ij} \in \R$ such that $(\sigma_0^*S_M)(t_i, e_j) = \frac{1}{p_i} Z_{ij} \mod \Z$. We define the homomorphism $\tilde{\tau}_T : T^{2k+1} \to \bH^{2k}(M)$ by extending $\tilde{\tau}_T(t_i) = \frac{1}{p_i}\sum_{k = 1}^b Z_{ik}\eta_k$ linearly. The property of the basis $\{ h_1, \ldots, h_b \}$ shows:
$$
\int_M \tilde{\tau}_T(t_i) \wedge \delta(\sigma_0(e_j)) = \frac{1}{p_i}Z_{ij}.
$$
Thus, the composition of $\tilde{\tau}_T$ with the projection $\bH^{2k}(M) \to \bH^{2k}(M)/\bH^{2k}(M)_\Z$ provides the homomorphism $\tau_T : T^{2k+1} \to \bH^{2k}(M)/\bH^{2k}(M)_\Z$, which satisfies
$$
(\sigma_0^*S_M)(t, e) = \int_M \tau_T(t) \wedge \delta(\sigma_0(e)) \mod \Z
$$
for $t \in T^{2k+1}$ and $e \in F^{2k+1}$. We can directly verify that $\tau_T$ is the unique homomorphism $T^{2k+1} \to \bH^{2k}(M)/\bH^{2k}(M)_\Z$ satisfying the condition above. If we choose a homomorphism $\tau_F : F^{2k+1} \to \bH^{2k}(M)/\bH^{2k}(M)_\Z$ and define $\tau : H^{2k+1}(M, \Z) \to \bH^{2k}(M)/\bH^{2k}(M)_\Z$ by $\tau = \tau_T \oplus \tau_F$, then $\sigma_0 - \tau$ is a harmonic splitting compatible with $\omega$.
\end{proof}

The following lemmas are properties of compatible harmonic splittings.

\begin{lem} \label{lem:splitting_basis_compatible}
Let $\omega : H^{2k+1}(M, \Z) \cong F^{2k+1} \oplus T^{2k+1}$ be a decomposition.

(a) For a basis $\{ e_1, \ldots, e_b \}$ of $F^{2k+1}$, there is a harmonic splitting $\sigma$ compatible with $\omega$ such that $\sigma^*S_M(e_i, e_j) = 0$ for $i \neq j$.

(b) There is a harmonic splitting $\sigma$ compatible with $\omega$ such that 
$$
2(\sigma^*S_M)(e, c) = 0
$$ 
for $e \in F^{2k+1}$ and $c \in H^{2k+1}(M, \Z)$.
\end{lem}

\begin{proof}
Note that (b) follows from (a), because $2 (\sigma^*S_M)(e, e) = 0$ for all $e \in F^{2k+1}$. To prove (a), we construct a basis of the vector space $\bH^{2k+1}(M)$ as follows. We take a harmonic splitting $\sigma_0$ compatible with $\omega$, and define the harmonic $(2k+1)$-forms $\epsilon_i$ by $\epsilon_i = \delta(\sigma_0(e_i))$. The commutativity of the diagram (\ref{formula:de_Rham_real_coefficients}) implies that $\{ \epsilon_1, \ldots, \epsilon_b \}$ forms a basis of the vector space $\bH^{2k+1}(M)$. 

We here consider the following linear map of the spaces of matrices:
$$
M(b, \R) \longrightarrow 
\{ A \in M(b, \R) |\ A + {}^t\!A = 0 \}, \quad
\Phi \mapsto \Phi G - {}^t(\Phi G),
$$
where $G = (G_{ij}) \in M(b, \R)$ is defined to be $G_{ij} = (\epsilon_i, \epsilon_j)_{L^2}$ by using the $L^2$-inner product $( \ , \ )_{L^2}$ on $A^{2k+1}(M)$. Because $( \ , \ )_{L^2}$ induces a positive definite inner product on $\bH^{2k+1}(M)$, the symmetric matrix $G$ is invertible. Hence the linear map is surjective.

Now, we take a skew-symmetric matrix $A = (A_{ij}) \in M(b, \R)$ such that 
$$
(\sigma_0^*S_M)(e_i, e_j) = A_{ij} \mod \Z
$$ 
for $i \neq j$. Then there exists a matrix $\Phi = (\Phi_{ij}) \in M(b, \R)$ such that $\Phi G - {}^t(\Phi G) = A$. By computations, we can verify
\begin{equation*}
\begin{split}
(\sigma_0^*S_M)(e_i, e_j) 
&= 
\sum_k \Phi_{ik}G_{kj} - \sum_k G_{ik}\Phi_{jk} \mod \Z \\
&=
\int_M
\left( \sum_{k} \Phi_{ik} * \delta(\sigma_0(e_k)) \right) \wedge 
\delta(\sigma_0(e_j)) \\
&\quad \quad 
-
\int_M
\left( \sum_{k} \Phi_{jk} * \delta(\sigma_0(e_k)) \right) \wedge 
\delta(\sigma_0(e_i)) \mod \Z.
\end{split}
\end{equation*}
Since $\delta(\sigma_0(e_k)) \in A^{2k+1}(M)$ is a harmonic form, so is $*\delta(\sigma_0(e_k)) \in A^{2k}(M)$. We define the homomorphism $\tilde{\tau}_F : F^{2k+1} \to \bH^{2k}(M)$ by the linear extension of $\tilde{\tau}_F(e_i) = \sum_k \Phi_{ik} * \delta(\sigma_0(e_k))$. Composing $\tilde{\tau}_F$ with the projection $\bH^{2k}(M) \to \bH^{2k}(M)/\bH^{2k}(M)_\Z$, we obtain $\tau_F : F^{2k+1} \to \bH^{2k}(M)/\bH^{2k}(M)_\Z$. For $i \neq j$ this homomorphism satisfies
$$
(\sigma_0^*S_M)(e_i, e_j) 
= 
\int_M \tau_F(e_i) \wedge \delta(\sigma_0(e_j))
- \int_M \tau_F(e_j) \wedge \delta(\sigma_0(e_i)) \mod \Z.
$$
We define $\tau : H^{2k+1}(M, \Z) \to \bH^{2k}(M)/\bH^{2k}(M)_\Z$ so that $\tau(e) = \tau_F(e)$ for $e \in F^{2k+1}$ and $\tau(t) = 0$ for $t \in T^{2k+1}$. Then the harmonic splitting $\sigma_0 - \tau$ has the property in (a), and is compatible with $\omega$ by Proposition \ref{prop:compatible_harmonic_splitting} (b).
\end{proof}

\begin{lem} \label{lem:decomposition_splitting}
Let $\omega : H^{2k+1}(M, \Z) \cong \Ker \ \bar{\omega} \oplus T^{2k+1}$ and $\omega' : H^{2k+1}(M, \Z) \cong \Ker \ \bar{\omega}' \oplus T^{2k+1}$ be decompositions of $H^{2k+1}(M, \Z)$. We put $\theta = \bar{\omega}' - \bar{\omega}$ and define a homomorphism $\theta : H^{2k+1}(M, \Z) \to T^{2k+1}$. Suppose that $\sigma$ and $\sigma'$ are harmonic splittings compatible with $\omega$ and $\omega'$, respectively. If we put $\tau = \sigma' - \sigma$, then we have:
\begin{align*}
L_M(t, \theta(e)) 
&= S_M(\tau(t), \sigma(e)), \\
({\sigma'}^*S_M)(e - \theta(e), \xi - \theta(\xi)) 
&=
(\sigma^*S_M)(e, \xi) - L_M(\theta(e), \theta(\xi)) \\
& \quad
+ S_M(\tau(e), \sigma(\xi)) - S_M(\tau(\xi), \sigma(e)),
\end{align*}
for $t \in T^{2k+1}$ and $e, \xi \in \Ker \ \bar{\omega}$.
\end{lem}

\begin{proof}
For convenience, we put $e' = e - \theta(e)$ for $e \in \Ker \ \bar{\omega}$. Clearly, $e'$ belongs to $\Ker \ \bar{\omega}'$. Because $\sigma$ and $\sigma'$ are harmonic, we have
$$
({\sigma'}^*S_M)(e', t) - (\sigma^*S_M)(e, t) 
= - L(\theta(e), t) + S_M(\sigma(e), \tau(t)).
$$
Since $\sigma$ and $\sigma'$ are compatible splittings, we have the first formula in this lemma. For the second formula, we also put $\xi' = \xi - \theta(\xi)$. We use the harmonicity of $\sigma$ and $\sigma'$ again to give
\begin{align*}
({\sigma'}^*S_M)(e', \xi') 
&= 
(\sigma^*S_M)(e, \xi) 
+ S_M(\tau(e), \sigma(\xi)) - S_M(\tau(\xi), \sigma(e)) \\
&
+ L_M(\theta(e), \theta(\xi)) 
- S_M(\tau(\theta(e)), \sigma(\xi)) + S_M(\tau(\theta(\xi)), \sigma(e)).
\end{align*}
Substitutions of the first formula complete the proof.
\end{proof}

\subsection{Decomposition of the smooth Deligne cohomology}

We introduce here some other notions related to $\G(M)$ and $S_M$.

\begin{dfn}
For a decomposition $\omega : H^{2k+1}(M, \Z) \cong F^{2k+1} \oplus T^{2k+1}$, we define $\G_\omega(M)$ to be the following subgroup in $\G(M)$:
$$
\G_\omega(M) 
= \Ker(\bar{\omega} \circ \chi)
= \{ f \in \G(M) |\ \chi(f) \in F^{2k+1} \}.
$$
\end{dfn}

We denote by $\til{\G}_\omega(M)$ the central extension of $\G_\omega(M)$ given by restricting $\til{\G}(M)$. Note that the linking form $L_M : T^{2k+1} \times T^{2k+1} \to \R/\Z$ gives a group 2-cocycle of $T^{2k+1}$. We also denote by $\til{T}^{2k+1}$ the central extension of $T^{2k+1}$ determined by $e^{2\pi\im L_M}$. 

Let $\til{\G}_\omega(M) \otimes \til{T}^{2k+1}$ be the quotient space of $\til{\G}_\omega(M) \times \til{T}^{2k+1}$ under the diagonal action of $\T$. The quotient space gives rise to a (topological) group by the multiplication $(\til{f}_1 \otimes \til{t}_1) \cdot (\til{f}_2 \otimes \til{t}_2) = (\til{f}_1 \til{f}_2) \otimes (\til{t}_1 \til{t}_2)$. In particular, the group is a central extension of $\G_\omega(M) \times T^{2k+1}$.

\begin{lem} \label{lem:decompose_central_ext}
For a decomposition $\omega : H^{2k+1}(M, \Z) \cong F^{2k+1} \oplus T^{2k+1}$, there are canonical isomorphisms $\G_\omega(M) \times T^{2k+1} \to \G(M)$ and $\til{\G}_\omega(M) \otimes \til{T}^{2k+1} \to \til{\G}(M)$ which make the following diagram commutative:
$$
\begin{CD}
1 @>>> 
\T @>>> 
\til{\G}_\omega(M) \otimes \til{T}^{2k+1} @>>> 
\G_\omega(M) \times T^{2k+1} @>>> 
1 \\
@. @| @VVV @VVV @. \\
1 @>>> 
\T @>>> 
\til{\G}(M) @>>>
\G(M) @>>>
1.
\end{CD}
$$
\end{lem}

\begin{proof}
We introduce the following maps:
\begin{align*}
I_\omega : \ \G_\omega(M) \times T^{2k+1} 
&\longrightarrow \G(M), &
(f, t) &\mapsto f + \sigma(t), \\
\til{I}_\omega : \ \til{\G}_\omega(M) \otimes \til{T}^{2k+1} 
&\longrightarrow \til{\G}(M), &
(f, u) \otimes (t, v) &\mapsto (f + \sigma(t), uv),
\end{align*}
where $\sigma$ is a harmonic splitting compatible with $\omega$. By Proposition \ref{prop:compatible_harmonic_splitting} (b), these maps are independent of the choice of $\sigma$. Clearly, $I_\omega$ is an isomorphism. We can see that $\til{I}_\omega$ is also an isomorphism, because $S_M(f_1 + \sigma(t_1), f_2 + \sigma(t_2)) = S_M(f_1, f_2) + L_M(t_1, t_2)$ holds for $f_i \in \G_\omega(M)$ and $t_i \in T^{2k+1}$.
\end{proof}

\begin{dfn} \label{dfn:Weyl_group_action}
We introduce the homomorphism
$$
s : H^{2k+1}(M, \Z) \longrightarrow 
\mathrm{Hom}(\bH^{2k}(M)/\bH^{2k}(M)_\Z, \R/\Z)
$$
by defining $s(\xi) : \bH^{2k}(M)/\bH^{2k}(M)_\Z \to \R/\Z$ to be 
$$
s(\xi)(\eta) = 
S_M(\eta, \sigma(\xi)) = 
\int_M \eta \wedge \delta(\sigma(\xi)) \mod \Z,
$$
where $\sigma$ is a harmonic splitting.
\end{dfn}

The homomorphism $s$ is independent of the choice of $\sigma$.

\begin{lem} \label{lem:Weyl_group_action}
(a) The homomorphism $s$ induces the isomorphism
$$
H^{2k+1}(M, \Z)/T^{2k+1} \cong \mathrm{Hom}(\bH^{2k}(M)/\bH^{2k}(M)_\Z, \R/\Z).
$$

(b) Let $\theta : H^{2k+1}(M, \Z) \to T^{2k+1}$ be a homomorphism such that $\theta(t) = 0$ for $t \in T^{2k+1}$. There is a unique homomorphism $\tau : T^{2k+1} \to \bH^{2k}(M)/\bH^{2k}(M)_\Z$ such that
$$
L_M(t, \theta(c)) = (s(c) \circ \tau)(t)
$$
for all $t \in T^{2k+1}$ and $c \in H^{2k+1}(M, \Z)$.
\end{lem}

\begin{proof}
We can prove (a), using the bases $\{ e_i \}$ and $\{ h_i \}$ in the proof of Proposition \ref{prop:compatible_harmonic_splitting}. Using the bases again, we can see the uniqueness of the homomorphism $\tau$ in (b). So, it is enough to show the existence. Let $\omega : H^{2k+1}(M, \Z) \cong F^{2k+1} \oplus T^{2k+1}$ be a decomposition. We have the other decomposition $\omega'$ by $\bar{\omega}' = \bar{\omega} + \theta$. Taking harmonic splittings $\sigma$ and $\sigma'$ compatible with $\omega$ and $\omega'$, respectively, we define $\tau_T : T^{2k+1} \to \bH^{2k}(M)/\bH^{2k}(M)_\Z$ by restricting the homomorphism $\sigma' - \sigma$. Lemma \ref{lem:decomposition_splitting} implies that $\tau_T$ has the property in (b).
\end{proof}


\section{Construction of admissible representations}
\label{sec:construction}

In this section, we construct projective unitary representations of the smooth Deligne cohomology group $\G(M)$, generalizing that of positive energy representations of $L\T$ given in \cite{P-S}. We also show that the representations are \textit{admissible} in the sense of Definition \ref{idfn:admissible}.

\bigskip

When $(\rho, H)$ is a projective representation of a group $G$, we denote by $(\til{\rho}, H)$ the corresponding linear representation of its central extension $\til{G}$, and vice verse. This identification will be used freely.

\subsection{Outline of construction}

We sketch the construction of representations of $\til{\G}(M)$. Suppose that the following data are given:
\begin{list}{}{\parsep=-2pt\topsep=4pt}
\item[(i)]
a decomposition 
$\omega : H^{2k+1}(M, \Z) \cong F^{2k+1} \oplus T^{2k+1}$;

\item[(ii)]
a homomorphism $\lambda : \bH^{2k}(M)/\bH^{2k}(M)_\Z \to \R/\Z$;

\item[(iii)]
a finite dimensional projective unitary representation $(\pi, V)$ of $T^{2k+1}$ with its cocycle $e^{2\pi\im L_M} : T^{2k+1} \times T^{2k+1} \to \T$.
\end{list}
We notice the following decomposition:
$$
\G(M) \cong 
(\bH^{2k}(M)/\bH^{2k}(M)_\Z) \times d^*(A^{2k+1}(M)) \times 
F^{2k+1} \times T^{2k+1}.
$$
Now, the construction proceeds as follows: First, we construct a representation $(\tilde{\rho}, H)$ of the central extension of $d^*(A^{2k+1}(M))$ as that of the \textit{Heisenberg group}. Second, we extend $(\til{\rho}, H)$ to the representation $(\til{\rho}_\lambda, H_\lambda)$ of the central extension of $\G^0(M) \cong (\bH^{2k}(M)/\bH^{2k}(M)_\Z) \times d^*(A^{2k+1}(M))$ by letting $\bH^{2k}(M)/\bH^{2k}(M)_\Z$ act on $H = H_\lambda$ through $\lambda$. Then, by using the representation of the subgroup $\til{\G}^0(M)$ in $\til{\G}_\omega(M)$, we obtain the induced representation $(\til{\rho}^\omega_\lambda, \H^\omega_\lambda)$ of $\til{\G}_\omega(M)$. Finally, we define the representation $(\til{\rho}^\omega_{\lambda, V}, \H^\omega_{\lambda, V})$ of $\til{\G}(M) \cong \til{\G}_\omega(M) \otimes \til{T}^{2k+1}$ to be the tensor product $\H^\omega_{\lambda, V} = \H^\omega_\lambda \otimes V$.

\medskip

If $M = S^1$, then the above decomposition of $\G(S^1) \cong L\T$ is
$$
\textstyle{
L \T \cong 
\T \times \{ \phi : S^1 \to \R |\ \int_{S^1} \phi(t) dt = 0 \} \times \Z.
}
$$
In this case, the construction of representations sketched above coincides with that of positive energy representations of $L\T$ presented in \cite{P-S, S}.

\subsection{The Heisenberg group}
\label{subsec:Heisenberg}

We recall a general construction of representations of Heisenberg groups \cite{P-S,S}.

\medskip

Let $V$ be a topological vector space over $\R$, and $S : V \times V \to \R$ a continuous skew-symmetric bilinear form which is non-degenerate in the following sense: any non-zero $v \in V$ has an element $v' \in V$ such that $S(v, v') \neq 0$. The \textit{Heisenberg group} associated to $(V, S)$ is defined to be the set $\tilde{V} = V \times \T$ endowed with the group multiplication
$$
(v, z) \cdot (v', z') = (v + v', zz' e^{\im S(v, v')}).
$$
By definition, $\tilde{V}$ is a central extension of $V$ by $\T$:
$$
\begin{CD}
1 @>>> \T @>>> \tilde{V} @>>> V @>>> 1.
\end{CD}
$$

\begin{prop}[\cite{P-S}] \label{prop:Heisenberg}
Let $\tilde{V}$ be the Heisenberg group associated to $(V, S)$. 

(a) Suppose that a continuous complex structure $J : V \to V$ obeys:
\begin{list}{}{\parsep=-2pt\topsep=4pt}
\item[(J1)]
$J$ is compatible with $S$, that is, $S(Jv, Jv') = S(v, v')$ for $v, v' \in V$;

\item[(J2)]
$J$ is positive, that is, $S(Jv, v) > 0$ for $v \in V$ such that $v \neq 0$.
\end{list}
Then there exists a unitary representation $(\tilde{\rho}, H)$ of $\tilde{V}$ on a Hilbert space $H$ such that the center $\T$ acts as the scalar multiplication. 

(b) The representation is continuous in the sense that $\tilde{\rho} : \tilde{V} \times H \to H$ is a continuous map.

(c) Suppose that $V$ is separable and complete with respect to the positive definite inner product $( \ , \ ) : V \times V \to \R$ defined by $(v, v') = S(Jv, v')$. Then the representation $(\tilde{\rho}, H)$ is irreducible.
\end{prop}

The unitary representation $(\tilde{\rho}, H)$ of $\til{V}$ corresponds to a projective unitary representation $(\rho, H)$ of $V$ with its cocycle $e^{\im S} : V \times V \to \T$. We briefly explain the construction of $(\rho, H)$ described in \cite{P-S}.

\smallskip

First, we describe the representation space $H$. Let $V_\C$ be the complexification $V_\C = V \otimes \C$. We extend the inner product on $V$ to the Hermitian inner product $( \ , \ ) : V_\C \times V_\C \to \C$ so as to be $\C$-linear in the first variable. (This convention differs from \cite{P-S,S}.) The $\C$-linear extension of $J$ provides us the decomposition $V_\C = W \oplus \overline{W}$, where $J$ acts on $W$ and $\overline{W}$ by $\im$ and $-\im$, respectively. To the symmetric algebra $S(W) = \bigoplus_{k \ge 0} S^k(W)$ of $W$, we introduce the Hermitian inner product $\langle \ , \ \rangle$ as follows:
$$
\langle w_1 \cdots w_k, w'_1 \cdots w'_k \rangle =
\sum_{\sigma \in \mathfrak{S}_k} 2^k 
\langle w_{\sigma(1)}, w'_1 \rangle \cdots
\langle w_{\sigma(k)}, w'_k \rangle,
$$
where $\mathfrak{S}_k$ denotes the symmetric group of degree $k$. Now, the representation space $H$ is defined to be the Hilbert space $\widehat{S}(W)$ obtained by completing $S(W)$.

Next, we consider a dense subspace in $H = \widehat{S}(W)$. For $\xi \in W$, the sequence $\{ \sum_{j = 0}^n \xi^j/j! \}_{n \in \N}$ in $S(W)$ converges to an element in $\widehat{S}(W)$, which we denote by $\epsilon_\xi$. Note that $\langle \epsilon_\xi, \epsilon_\eta \rangle = e^{\langle \xi, \eta \rangle}$. Note also that $\epsilon_{\xi_1}, \ldots, \epsilon_{\xi_n} \in H$ are linearly independent, provided that $\xi_1, \ldots, \xi_n \in W$ are distinct. Let $E$ be the subspace in $\widehat{S}(W)$ generated by $\{ \epsilon_\xi |\ \xi \in W \}$. The key fact is that the closure $\bar{E}$ of $E$ coincides with $\widehat{S}(W)$.

Finally, we describe the action of $V$ on $H = \widehat{S}(W)$. For $v_+ \in W$ and $v_- \in \overline{W}$, we define the linear map $\rho(v_+ + v_-) : \ E \to E$ by
$$
\rho(v_+ + v_-) \epsilon_\xi = 
\exp\left(
- \frac{1}{2} \langle v_+, \overline{(v_-)} \rangle 
- \langle \xi, \overline{(v_-)} \rangle \right)
\epsilon_{\xi + v_+}.
$$
We can verify $\rho(v)\rho(v') \epsilon_\xi = e^{\im (v, J \bar{v}')} \rho(v + v') \epsilon_\xi$ for $v, v' \in V_\C$, so that we have a projective representation $\rho : V_\C \times E \to E$. In general, the linear map $\rho(v) : E \to E$ is unbounded. However, if $v$ belongs to $V \subset V_\C$, then $\rho(v) : E \to E$ is isometric. Consequently, the extension $\rho(v) : H \to H$ gives the projective unitary representation $\rho : V \times H \to H$. For $\xi \in W$ the map $\rho(\cdot)\epsilon_\xi : V_\C \to E$ is continuous. Thus, $\rho : V \times H \to H$ gives rise to a continuous map.

\subsection{Construction}

\begin{lem} \label{lem:Heisenberg}
We can construct an irreducible continuous projective unitary representation $(\rho, H)$ of $d^*(A^{2k+1}(M))$ on a Hilbert space $H$ with its cocycle $e^{2\pi\im S_M}$. 
\end{lem}

\begin{proof}
Note that the group 2-cocycle 
$$
S_M : d^*(A^{2k+1}(M)) \times d^*(A^{2k+1}(M)) \to \R/\Z
$$
is expressed as $S_M(\nu, \nu') = \int_M \nu \wedge d \nu' \mod \Z$. Combining Proposition \ref{prop:topology_differential_forms} with Proposition \ref{prop:Heisenberg}, we obtain a continuous projective unitary representation $(\rho, V)$ of the completion $V$ of $d^*(A^{2k+1}(M))$. By restricting the action of $V$, we obtain the continuous projective unitary representation of $d^*(A^{2k+1}(M))$ stated in this lemma. The representation is irreducible, since $\rho$ is continuous and $d^*(A^{2k+1}(M))$ is dense in $V$.
\end{proof}

As in Section \ref{sec:introduction}, we put $\mathcal{X}(M) = \Hom(\bH^{2k}(M)/\bH^{2k}(M)_\Z, \R/\Z)$.

\begin{lem} \label{lem:representation_identity_component}
For $\lambda \in \mathcal{X}(M)$, we can construct an irreducible continuous projective unitary representation $(\rho_\lambda, H_\lambda)$ of $\G^0(M) = A^{2k}(M)/A^{2k}(M)_\Z$ on a Hilbert space $H_\lambda$ with its cocycle $e^{2\pi\im S_M}$.
\end{lem}

This establishes Lemma \ref{ilem:identity_component} in Section \ref{sec:introduction}.

\begin{proof}
We put $H_\lambda = H$ and define $\rho_\lambda : \G^0(M) \times H_\lambda \to H_\lambda$ by $\rho_\lambda(\eta, \nu) = \lambda(\eta) \rho(\nu)$ for $(\eta, \nu) \in \bH^{2k}(M)/\bH^{2k}(M)_\Z \times d^*(A^{2k+1}(M)) \cong \G^0(M)$. Then $(\rho_\lambda, H_\lambda)$ gives rise to a continuous projective unitary representation with its cocycle $e^{2\pi\im S_M}$, since the group 2-cocycle $S_M$ on $\G^0(M)$ has the expression
$$
S_M((\eta, \nu), (\eta', \nu')) = \int_M \nu \wedge d\nu' \mod \Z.
$$
Clearly, the restriction of $(\rho_\lambda, H_\lambda)$ to the subgroup $d^*(A^{2k+1}(M)) \subset \G^0(M)$ coincides with $(\rho, H)$. Hence $(\rho_\lambda, H_\lambda)$ is irreducible.
\end{proof}

\begin{dfn} 
We fix a decomposition $\omega : H^{2k+1}(M, \Z) \cong F^{2k+1} \oplus T^{2k+1}$. For $\lambda \in \mathcal{X}(M)$, we define $(\til{\rho}^\omega_\lambda, \H^\omega_\lambda)$ to be the unitary representation of $\til{\G}_\omega(M)$ induced from the unitary representation $(\til{\rho}_\lambda, H_\lambda)$ of the subgroup $\til{\G}^0(M)$.
\end{dfn}

We give here a concrete description of the projective unitary representation $(\rho^\omega_\lambda, \H^\omega_\lambda)$. The unitary action of $\til{\G}^0(M)$ on $H_\lambda$ defines the associated Hermitian vector bundle $\til{\G}^\omega(M) \times_{\til{\rho}^\omega_\lambda} H_\lambda$ over $\til{\G}^\omega(M)/\til{\G}^0(M)$ equivariant under $\til{\G}^\omega(M)$. By definition, the representation space $\H^\omega_\lambda$ consists of the square summarable sections of the vector bundle. Note that $\chi : \G_\omega(M) \to F^{2k+1}$ induces an isomorphism of groups $\til{\G}_\omega(M)/\til{\G}^0(M) \cong F^{2k+1}$. A choice of a splitting $\sigma : H^{2k+1}(M, \Z) \to \G(M)$ trivializes the vector bundle as follows:
$$
F^{2k+1} \times H_\lambda \to
\til{\G}_\omega(M) \times_{\til{\rho}^\omega_\lambda} H_\lambda, \quad
(\xi, v) \mapsto [(\sigma(\xi), 1), v].
$$
Then the representation space $\H^\omega_\lambda$ is the Hilbert space direct sum
\begin{equation} \label{formula:concrete_representation_space}
\H^\omega_\lambda = \bigg\{
\Phi : F^{2k+1} \to H \bigg|\
\sum_{\xi \in F^{2k+1}} \p{\Phi(\xi)}^2 < + \infty
\bigg\},
\end{equation}
and the action of $f \in \G_\omega(M)$ on $\Phi \in \H^\omega_\lambda$ is
$$
(\rho^\omega_\lambda(f)\Phi)(\xi)
= \\
e^{2\pi\im \left\{
S_M(\sigma(c), \sigma(\xi) + \alpha) 
- S_M(\sigma(c), \sigma(c))
+ 2S_M(\alpha, \sigma(\xi))
\right\}}
\rho_{\lambda}(\alpha) \Phi(\xi - c),
$$
where $c = \chi(f) \in F^{2k+1}$ and $\alpha = f - \sigma(\chi(f)) \in \G^0(M)$. If the splitting $\sigma$ is harmonic, then we have a simpler expression using the homomorphism $s$ in Definition \ref{dfn:Weyl_group_action}:
\begin{equation}
(\rho^\omega_\lambda(f)\Phi)(\xi)
=
e^{2\pi\im \left\{
S_M(\sigma(c), \sigma(\xi) + \alpha) 
- S_M(\sigma(c), \sigma(c))
\right\}}
\rho_{\lambda + 2 s(\xi)}(\alpha) \Phi(\xi - c).
\label{formula:concrete_representation_action}
\end{equation}

\begin{lem} \label{lem:branching_rule}
The projective unitary representation of $\G^0(M)$ given by the restriction of $(\rho^\omega_\lambda, \H^\omega_\lambda)$ is the Hilbert space direct sum $\widehat{\oplus}_{\xi \in F^{2k+1}} H_{\lambda + 2s(\xi)}$.
\end{lem}

\begin{proof}
We have $(\rho^\omega_\lambda(\alpha)\Phi)(\xi) = \rho_{\lambda + 2 s(\xi)}(\alpha)\Phi(\xi)$ for $\alpha \in \G^0(M)$ by (\ref{formula:concrete_representation_action}).
\end{proof}

\medskip

Now we construct projective representations of $\G(M)$. 

\begin{dfn} \label{dfn:representation}
We suppose that the following data are given:
\begin{list}{}{\parsep=-2pt\topsep=4pt}
\item[(i)]
a decomposition 
$\omega : H^{2k+1}(M, \Z) \cong F^{2k+1} \oplus T^{2k+1}$;

\item[(ii)]
a homomorphism $\lambda \in \mathcal{X}(M)$;

\item[(iii)]
a finite dimensional projective unitary representation $(\pi, V)$ of $T^{2k+1}$ with its cocycle $e^{2\pi\im L_M}$.
\end{list}
We define the projective unitary representation $(\rho^\omega_{\lambda, V}, \H^\omega_{\lambda, V})$ of $\G(M)$ with its cocycle $e^{2\pi\im S_M}$ as follows: The Hilbert space  $\H^\omega_{\lambda, V}$ is the algebraic tensor product $\H^\omega_{\lambda, V} = \H^\omega_\lambda \otimes V$. The action $\rho^\omega_{\lambda, V} : \G(M) \times \H^\omega_{\lambda, V} \to \H^\omega_{\lambda, V}$ is 
$$
\rho^\omega_{\lambda, V}(f) (\Phi \otimes v) 
= 
(\rho^\omega_\lambda(f - \sigma(t)) \Phi) \otimes (\pi(t) v),
$$
where $\sigma$ is a harmonic splitting compatible with $\omega$, and $t = \bar{\omega} \circ \chi(f) \in T^{2k+1}$.
\end{dfn}

In the definition above, the isomorphism $\til{\G}(M) \cong \til{\G}_\omega(M) \otimes \til{T}^{2k+1}$ in Lemma \ref{lem:decompose_central_ext} is used. Hence $\rho^\omega_{\lambda, V}$ is independent of the choice of $\sigma$. We can also verify the independence directly by means of Proposition \ref{prop:compatible_harmonic_splitting} (b).

\begin{rem}
Since $V$ is finite dimensional, the algebraic tensor product $\H^\omega_\lambda \otimes V$ is naturally isomorphic to the Hilbert space tensor product $\H^\omega_\lambda \widehat{\otimes} V$, which is the completion of $\H^\omega_\lambda \otimes V$ with respect to the unique inner product such that $\langle \phi_1 \otimes v_1, \phi_2 \otimes v_2 \rangle = \langle \phi_1, \phi_2 \rangle \langle v_1, v_2 \rangle$ for elementary tensors.
\end{rem}

\begin{prop} \label{prop:branching_rule}
The projective unitary representation $\H^\omega_{\lambda, V}|_{\G^0(M)}$ of $\G^0(M)$ given by restricting $(\rho^\omega_{\lambda, V}, \H^\omega_{\lambda, V})$ is
$$
\H^\omega_{\lambda, V}|_{\G^0(M)} \cong
\widehat{\bigoplus}_{\xi \in F^{2k+1}} 
\bigoplus^{\mathrm{dim} V} 
H_{\lambda + 2s(\xi)}.
$$
Therefore $(\rho^\omega_{\lambda, V}, \H^\omega_{\lambda, V})$ is admissible in the sense of Definition \ref{idfn:admissible}.
\end{prop}

\begin{proof}
The (unitary) equivalence is clear by Lemma \ref{lem:branching_rule}. Now we can readily see that $(\rho^\omega_{\lambda, V}, \H^\omega_{\lambda, V})$ is admissible, using Lemma \ref{lem:Weyl_group_action} (a).
\end{proof}

\section{Classification of admissible representations}
\label{sec:classification}

In this section, we prove Theorem \ref{ithm:classification} in Section \ref{sec:introduction} --- the classification of admissible representations of $\G(M)$. For this aim, we first classify the admissible representations $(\rho^\omega_{\lambda, V}, \H^\omega_{\lambda, V})$ in Definition \ref{dfn:representation}. We also give a condition for them to be irreducible. Then we prove that any admissible representation is equivalent to the direct sum of some $(\rho^\omega_{\lambda, V}, \H^\omega_{\lambda, V})$, by modifying the given equivalence of projective representations of $\G^0(M)$.

\subsection{Dependence on free part}

First of all, we investigate how $(\rho^\omega_{\lambda, V}, \H^\omega_{\lambda, V})$ depends on the choice of $\omega$.

\begin{dfn} \label{dfn:twisting_projective_representation}
Let $(\pi, V)$ be a projective unitary representation of $T^{2k+1}$ with its cocycle $e^{2\pi\im L_M}$. 

(a) For a homomorphism $\mu : T^{2k+1} \to \R/\Z$, we define the projective unitary representation $(\pi_\mu, V_\mu)$ of $T^{2k+1}$ with its cocycle $e^{2\pi\im L_M}$ by setting 
\begin{align*}
V_\mu &= V, &
\pi_\mu(t) &= e^{2\pi\im \mu(t)} \pi(t).
\end{align*} 

(b) For homomorphisms $\lambda \in \mathcal{X}(M)$ and $\theta : H^{2k+1}(M, \Z) \to T^{2k+1}$ such that $\theta(t) = 0$ for $t \in T^{2k+1}$, we define the homomorphism $\mu : T^{2k+1} \to \R/\Z$ by $\mu = \lambda \circ \tau$, where $\tau : T^{2k+1} \to \bH^{2k}(M)/\bH^{2k}(M)_\Z$ is the homomorphism in Lemma \ref{lem:Weyl_group_action} (b). We write $(\pi_{(\lambda, \theta)}, V_{(\lambda, \theta)})$ for $(\pi_\mu, V_\mu)$.
\end{dfn}

\begin{prop} \label{prop:decomposition_dependence}
If $\theta: H^{2k+1}(M, \Z) \to T^{2k+1}$ is a homomorphism such that $\theta(t) = 0$ for $t \in T^{2k+1}$, then $(\rho^\omega_{\lambda, V}, \H^\omega_{\lambda, V})$ and $(\rho^{\omega + \theta}_{\lambda, V_{(\lambda, \theta)}}, \H^{\omega + \theta}_{\lambda, V_{(\lambda, \theta)}})$ are unitary equivalent.
\end{prop}

\begin{proof}
To suppress notations, we put $\bar{\omega}' = \bar{\omega} + \theta$ and $(\pi', V') = (\pi_{(\lambda, \theta)}, V_{(\lambda, \theta)})$. A choice of harmonic splittings $\sigma$ and $\sigma'$ compatible with $\omega$ and $\omega'$ respectively provides us the following identifications of the representation spaces:
\begin{align*}
\H^\omega_{\lambda, V} &=
\big\{
\Phi : \Ker \ \bar{\omega} \to H_\lambda \otimes V \big|\
\mbox{$\sum_{\xi \in \Ker \ \bar{\omega}} \p{\Phi(\xi)}^2 < + \infty$}
\big\}, \\
\H^{\omega'}_{\lambda, V'} &=
\big\{
\Phi' : \Ker \ \bar{\omega}' \to H_\lambda \otimes V \big|\
\mbox{$\sum_{\xi' \in \Ker \ \bar{\omega}'} \p{\Phi'(\xi')}^2 < + \infty$}
\big\}.
\end{align*}
We put $\tau = \sigma' - \sigma$. For $\Phi \in \H^\omega_{\lambda, V}$, we define $\mathcal{F}\Phi \in \H^{\omega'}_{\lambda, V'}$ to be
$$
\left(\mathcal{F} \Phi \right) (\xi') = 
e^{2\pi\im r(\xi)} 
\left(\id_{H_\lambda} \otimes \pi(\theta(\xi) \right) \Phi(\xi),
$$
where $\xi \in \Ker \ \bar{\omega}$ is given by $\xi = \xi'+\theta(\xi')$ and 
$$
r(\xi) = - \lambda \circ \tau(\xi) - S_M(\tau(\xi), \sigma(\xi)) + \lambda \circ \tau(\theta(\xi)).
$$
Using Lemma \ref{lem:decomposition_splitting}, we can check that $\mathcal{F} : \H^\omega_{\lambda, V} \to \H^{\omega'}_{\lambda, V'}$ is a unitary isomorphism satisfying $\rho^{\omega'}_{\lambda, V'} (f) \circ \mathcal{F} = \mathcal{F} \circ \rho^\omega_{\lambda, V}(f)$ for all $f \in \G(M)$.
\end{proof}

The proposition above allows us to fix a decomposition $\omega$ in studying the projective representations $(\rho^\omega_{\lambda, V}, \H^\omega_{\lambda, V})$.

\subsection{The space of intertwiners}

When (pre-)Hilbert spaces $H$ and $H'$ are given, we denote by $\Hom(H, H')$ the abstract linear space of continuous (bounded) linear maps. When $H$ and $H'$ are the representation spaces of some representations of a group $G$, we denote by $\Hom_G(H, H')$ the subspace consisting of $G$-equivariant maps.

\medskip

We begin with a simple property of a Hilbert space direct sum. 

\begin{lem} \label{lem:direct_sum_injection}
Let $\{ H_i \}_{i \in \N}$ be a sequence of Hilbert spaces. For another Hilbert space $H$, the following linear map given by the restriction is injective:
$$
\Hom(\widehat{\oplus}_{i \in \N} H_i, H) \to \prod_{i \in \N} \Hom(H_i, H).
$$
\end{lem}

\begin{proof}
The linear map factors as follows:
$$
\Hom(\widehat{\oplus}_i H_i, H) \to 
\Hom(\oplus_i H_i, H) \to
\prod_i \Hom(H_i, H).
$$
The first map is bijective, because $\oplus_i H_i$ is dense in $\widehat{\oplus}_i H_i$. The second map is obviously injective.
\end{proof}

\begin{lem} \label{lem:Frobenius_reciprocity}
There is a monomorphism
$$
\varpi: \
\Hom_{\til{\G}_\omega(M)}(\H^\omega_\lambda, \H^\omega_{\lambda'}) \to 
\Hom_{\til{\G}^0(M)}(\widehat{\oplus}_{\xi \in F^{2k+1}} H_{\lambda + 2 s(\xi)}, H_{\lambda'}).
$$
\end{lem}

\begin{proof}
To begin with, we define the linear map $\varpi$ as follows. Recall Lemma \ref{lem:branching_rule}. For $\mathcal{F} : \H^\omega_\lambda \to \H^\omega_{\lambda'}$ and $\Phi \in \H^\omega_\lambda$ we express $\mathcal{F}(\Phi) \in H^\omega_{\lambda'} = \widehat{\oplus}_\xi H_{\lambda' + 2s(\xi)}$ as $\mathcal{F}(\Phi) = \{ \mathcal{F}(\Phi)(\xi) \}_{\xi \in F^{2k+1}}$. We define $\varpi\mathcal{F}(\Phi) \in H_{\lambda'}$ to be $\varpi\mathcal{F}(\Phi) = \mathcal{F}(\Phi)(0)$. We can readily see that the linear map $\varpi\mathcal{F} : \widehat{\oplus}_{\xi} H_{\lambda + 2 s(\xi)} \to H_{\lambda'}$ is continuous. If $\mathcal{F}$ commutes with the action of $\G^\omega(M)$, then $\varpi\mathcal{F}$ commutes with the action of $\G^0(M)$. Hence the assignment $\mathcal{F} \mapsto \varpi\mathcal{F}$ gives rise to the $\varpi$.

Suppose that $\mathcal{F} \in \Hom_{\til{\G}_\omega(M)}(\H^\omega_\lambda, \H^\omega_{\lambda'})$ is such that $\varpi\mathcal{F} = 0$. We here take a compatible harmonic splitting $\sigma$ and apply (\ref{formula:concrete_representation_action}). Then we have
$$
0 = \varpi\mathcal{F}(\rho^\omega_\lambda(\sigma(c))\Phi) =
\left(\rho^\omega_{\lambda'}(\sigma(c))\mathcal{F}(\Phi)\right)(0) =
e^{-2\pi\im\sigma^*S_M(c, c)} \mathcal{F}(\Phi)(-c)
$$
for any $\Phi \in \H^\omega_\lambda$ and $c \in F^{2k+1}$. Thus $\mathcal{F} = 0$, so that $\varpi$ is injective.
\end{proof}

\begin{lem} \label{lem:distinction_representations_indentity_component}
If $\lambda \neq \lambda'$, then $\Hom_{\til{\G}^0(M)}(H_\lambda, H_{\lambda'}) = 0$.
\end{lem}

\begin{proof}
Let $F : H_\lambda \to H_{\lambda'}$ be a continuous linear map equivariant under the action of $\til{\G}^0(M)$. For $\eta \in \bH^{2k}(M)/\bH^{2k}(M)_\Z$ and $v \in H_\lambda$ we have 
$$
0 = F(\til{\rho}_\lambda(\eta, 1)v) - \til{\rho}_{\lambda'}(\eta, 1)F(v) 
= ( e^{2\pi\im\lambda(\eta)} - e^{2\pi\im\lambda'(\eta)} )F(v).
$$
Therefore $\lambda \neq \lambda'$ implies $F = 0$.
\end{proof}

By this lemma, $(\rho_\lambda, H_\lambda)$ and $(\rho_{\lambda'}, H_{\lambda'})$ are inequivalent for $\lambda \neq \lambda'$.

\begin{prop} \label{prop:intertwiners_I}
We suppose that $\lambda' \neq \lambda + 2s(\xi)$ for any $\xi \in F^{2k+1}$.

(a) $\Hom_{\til{\G}_\omega(M)}(\H^\omega_\lambda, \H^\omega_{\lambda'}) = 0$.

(b) $\Hom_{\til{\G}(M)}(\H^\omega_{\lambda, V}, \H^\omega_{\lambda', V'}) = 0$ for all $(\pi, V)$ and $(\pi', V')$.
\end{prop}

\begin{proof}
By Lemma \ref{lem:direct_sum_injection}, \ref{lem:Frobenius_reciprocity} and \ref{lem:distinction_representations_indentity_component}, we obtain (a). If the dimensions of $V$ and $V'$ are $d$ and $d'$, respectively, then the restriction of actions induces an injection:
$$
\Hom_{\til{\G}(M)}(\H^\omega_{\lambda, V}, \H^\omega_{\lambda', V'}) 
\longrightarrow
\bigoplus^{dd'} 
\Hom_{\til{\G}_\omega(M)}(\H^\omega_\lambda, \H^\omega_{\lambda'}).
$$
Hence (a) implies (b).
\end{proof}

\subsection{Equivalence}

The following proposition (``Schur's lemma'') plays a pivotal role in the subsections below. The proof is postponed in Appendix.

\begin{prop} \label{prop:Schur_lemma}
A unitary representation $(\rho, H)$ of a group $G$ on a Hilbert space $H$ is irreducible if and only if $\End_G(H) = \C \cdot \id_H$.
\end{prop}

As is seen in Lemma \ref{lem:representation_identity_component}, the unitary representation $(\til{\rho}_\lambda, H_\lambda)$ of $\til{\G}^0(M)$ is irreducible. Hence $\End_{\til{\G}^0(M)}(H_\lambda) \cong \C$.

\begin{lem} \label{lem:irreducible_free}
$\End_{\til{\G}_\omega(M)}(\H^\omega_\lambda) = \C \cdot \id_{\H^\omega_\lambda}$.
\end{lem}

\begin{proof}
Apparently, $\End_{\til{\G}_\omega(M)}(\H^\omega_\lambda) \supset \C \cdot \id_{\H^\omega_\lambda}$. On the other hand, Lemma \ref{lem:direct_sum_injection}, \ref{lem:Frobenius_reciprocity} and \ref{lem:distinction_representations_indentity_component} yield the injection $\End_{\til{\G}_\omega(M)}(\H^\omega_\lambda) \to \End_{\til{\G}^0(M)}(H_\lambda) \cong \C$.
\end{proof}

\begin{prop} \label{prop:intertwiners_II}
$\Hom_{\til{\G}(M)}(\H^\omega_{\lambda, V}, \H^\omega_{\lambda, V'}) \cong 
\Hom_{\til{T}^{2k+1}}(V, V')$.
\end{prop}

\begin{proof}
The restriction of actions induces the obvious injection
$$
\Hom_{\til{\G}(M)}(\H^\omega_{\lambda, V}, \H^\omega_{\lambda, V'}) 
\to 
\Hom_{\til{\G}_\omega(M)}
(\H^\omega_{\lambda, V}|_{\til{\G}_\omega(M)}, 
\H^\omega_{\lambda, V'}|_{\til{\G}_\omega(M)}).
$$
We have $\H^\omega_{\lambda, V}|_{\til{\G}_\omega(M)} = \H^\omega_\lambda \otimes V = \oplus^{\mathrm{dim} V} \H^\omega_\lambda$. So Lemma \ref{lem:irreducible_free} leads to
$$
\Hom_{\til{\G}_\omega(M)}
(\H^\omega_{\lambda, V}|_{\til{\G}_\omega(M)}, 
\H^\omega_{\lambda, V'}|_{\til{\G}_\omega(M)}) 
\cong
(\C \cdot \id_{\H^\omega_\lambda}) \otimes \Hom(V, V').
$$
We now express $\mathcal{F} \in \Hom_{\til{\G}_\omega(M)}(\H^\omega_{\lambda, V}|_{\til{\G}_\omega(M)}, \H^\omega_{\lambda, V'}|_{\til{\G}_\omega(M)})$ as $\mathcal{F} = \id_{\H^\omega_\lambda} \otimes F$, where $F : V \to V'$ is a linear map. It is direct to see that: $\mathcal{F}$ belongs to $\Hom_{\til{\G}_\omega(M)}(\H^\omega_{\lambda, V}, \H^\omega_{\lambda, V'})$ if and only if $F$ belongs to $\Hom_{\til{T}^{2k+1}}(V, V')$.
\end{proof}

\begin{thm}
There exists an equivalence between $\H^\omega_{\lambda, V}$ and $\H^\omega_{\lambda', V'}$ if and only if:
\begin{list}{}{\parsep=-2pt\topsep=4pt}
\item[(i)]
there is an equivalence between $(\pi, V)$ and $(\pi', V')$; and

\item[(ii)]
there is $\xi \in F^{2k+1}$ such that $\lambda ' = \lambda + 2s(\xi)$.
\end{list}
\end{thm}

\begin{proof}
To prove the ``if'' part, we take $\xi \in F^{2k+1}$ and a compatible harmonic splitting $\sigma$. From $(\rho^\omega_{\lambda, V}, \H^\omega_{\lambda, V})$, we construct the other projective representation $(\rho', \H')$ by setting $\H' = \H^\omega_{\lambda, V}$ and $\rho'(f) = \rho^\omega_{\lambda, V}(\sigma(\xi))^{-1} \circ \rho^\omega_{\lambda, V}(f) \circ \rho^\omega_{\lambda,V}(\sigma(\xi))$. Apparently, $(\rho', \H')$ is unitary equivalent to $(\rho^\omega_{\lambda, V}, \H^\omega_{\lambda, V})$. We obtain 
$$
\rho'(f) 
= e^{4\pi\im S_M(f, \sigma(\xi))} 
\rho^\omega_{\lambda, V}(f)
= e^{4\pi\im (\sigma^*S_M)(\chi(f), \xi)} 
\rho^\omega_{\lambda + 2s(\xi), V}(f)
$$
by the help of (\ref{formula:concrete_representation_action}). We can suppose that $\sigma$ is taken such as in Lemma \ref{lem:splitting_basis_compatible} (b). Then $\rho' = \rho^\omega_{\lambda + 2s(\xi), V}$, so that $\H^\omega_{\lambda, V}$ and $\H^\omega_{\lambda + 2s(\xi), V}$ are unitary equivalent. Now the ``only if'' part follows from Proposition \ref{prop:intertwiners_I} and \ref{prop:intertwiners_II}. 
\end{proof}

\subsection{Irreducibility}

\begin{thm} \label{thm:irreducibility}
$(\rho^\omega_{\lambda, V}, \H^\omega_{\lambda, V})$ is irreducible if and only if $(\pi, V)$ is irreducible.
\end{thm}

\begin{proof}
The theorem directly follows from Proposition \ref{prop:Schur_lemma} and \ref{prop:intertwiners_II}.
\end{proof}

We recall here a fact in the theory of projective representations of finite groups \cite{Kar}. Let $\Gamma$ be a finite group, and $\alpha : \Gamma \times \Gamma \to \T$ a group 2-cocycle. An element $\gamma \in \Gamma$ is said to be \textit{$\alpha$-regular} if $\alpha(\gamma, \eta) = \alpha(\eta, \gamma)$ for all $\eta \in \Gamma$ such that $\gamma \eta = \eta \gamma$. It is known that the number of the equivalence classes of finite dimensional irreducible projective representations of $\Gamma$ with their cocycle $\alpha$ coincides with that of the $\alpha$-regular elements.

\begin{cor} \label{cor:irreducibility}
If a decomposition $\omega : H^{2k+1}(M, \Z) \cong F^{2k+1} \oplus T^{2k+1}$ is fixed, then the number of the equivalence classes of irreducible representations $(\rho^\omega_{\lambda, V}, \H^\omega_{\lambda, V})$ in Definition \ref{dfn:representation} is $2^br$, where $b = b_{2k+1}(M)$ is the Betti number and $r$ is the number of the elements of the set $\{ t \in H^{2k+1}(M, \Z) |\ 2 t = 0 \}$.
\end{cor}

\begin{proof}
An element $t \in T^{2k+1}$ is $e^{2\pi\im L_M}$-regular if and only if $L_M(2t, c) = 0$ for all $c \in T^{2k+1}$. The \Poincare duality implies that: if $L_M(t', c) = 0$ for all $c \in T^{2k+1}$, then we have $t' = 0$. Hence the number of the equivalence classes of irreducible projective representations of $T^{2k+1}$ with their cocycle $e^{2\pi\im L_M}$ is $r$. Now, by Lemma \ref{lem:Weyl_group_action} (a), we can identify the cokernel of $2s : F^{2k+1} \to \mathcal{X}(M)$ with $(\Z/2\Z)^b$, which completes the proof.
\end{proof}

\subsection{Classification}

In this subsection we prove Theorem \ref{ithm:classification} in Section \ref{sec:introduction}. 

To begin with, we give a characterization of representations of $\til{\G}_\omega(M)$.

\begin{lem} \label{lem:characterize_rep_free_part}
Let $\omega : H^{2k+1}(M, \Z) \cong F^{2k+1} \oplus T^{2k+1}$ be a decomposition. A projective representation $\rho : \G_\omega(M) \times \H \to \H$ of $\G_\omega(M)$ with its cocycle $e^{2\pi\im S_M}$ is, as data, equivalent to the following projective representations:
\begin{align}
\rho_0 &: \G^0(M) \times \H \to \H, &
\rho_0(\alpha) \circ \rho_0(\alpha') &= 
e^{2\pi\im S_M(\alpha, \alpha')} \rho_0(\alpha + \alpha') \\
\rho_\Z &: F^{2k+1} \times \H \to \H, &
\rho_\Z(e) \circ \rho_\Z(e') &= 
e^{2\pi\im \sigma^*S_M(e, e')} \rho_\Z(e + e')
\label{formula:characterize_rep_free_part_II}
\end{align}
satisfying the compatibility condition:
\begin{equation}
\rho_0(\alpha) \circ \rho_\Z(e) = 
e^{4\pi\im S_M(\alpha, \sigma(e))} \rho_\Z(e) \circ \rho_0(\alpha),
\label{formula:characterize_rep_free_part_III}
\end{equation}
where $\sigma$ is a fixed harmonic splitting compatible with $\omega$.
\end{lem}

\begin{proof}
Recall the isomorphism $\G_\omega(M) \cong \G^0(M) \times F^{2k+1}$ induced by $\sigma$. Suppose that a projective representation $(\rho, \H)$ of $\G_\omega(M)$ is given. We have the projective representation $(\rho_0, \H)$ of $\G^0(M)$ by the restriction. We also have the projective representation $(\rho_\Z, \H)$ of $F^{2k+1}$ by setting $\rho_\Z(e) = \rho(\sigma(e))$. Then (\ref{formula:characterize_rep_free_part_III}) is clear. Conversely, if we put $\rho(\alpha + \sigma(e)) = e^{2\pi\im S_M(\alpha, \sigma(e))} \rho_0(\alpha) \circ \rho_\Z(e)$, then we recover $(\rho, \H)$ form $\rho_0$ and $\rho_\Z$.
\end{proof}

Now, we start the proof of Theorem \ref{ithm:classification}. Let $(\rho, \H)$ be an admissible representation. By definition, there is an equivalence of representations of $\G^0(M)$:
$$
\theta : \
\widehat{\bigoplus}_{\lambda \in \mathcal{X}(M)} 
\mathcal{V}(\lambda)
\longrightarrow
\H|_{\G^0(M)},
$$
where we put $\mathcal{V}(\lambda) = \oplus^{m(\lambda)} H_\lambda = H_\lambda \otimes \C^{m(\lambda)}$ to suppress notations. By virtue of Proposition \ref{appendix_prop:unitary_equivalence}, we can assume that $\theta$ is a \textit{unitary} equivalence. This unitary map induces the projective unitary representation
$$
\rho^\theta : \
\G(M) \times
\widehat{\bigoplus}_{\lambda \in \mathcal{X}(M)} 
\mathcal{V}(\lambda)
\longrightarrow
\widehat{\bigoplus}_{\lambda \in \mathcal{X}(M)} 
\mathcal{V}(\lambda)
$$
by $\rho^\theta(f) = \theta^{-1} \circ \rho(f) \circ \theta$. We take and fix a decomposition $\omega : H^{2k+1}(M, \Z) \cong F^{2k+1} \oplus T^{2k+1}$ and a compatible harmonic splitting $\sigma : H^{2k+1}(M, \Z) \to \G(M)$.

\begin{lem} \label{lem:Weyl_shift}
For $e \in F^{2k+1}$, $\rho^\theta(\sigma(e))$ restricts to a unitary isomorphism
$$
\rho^\theta(\sigma(e)) : \
\mathcal{V}(\lambda)
\longrightarrow
\mathcal{V}(\lambda + 2s(e)).
$$
Thus, $m(\lambda) = m(\lambda + 2s(e))$ for $\lambda \in \mathcal{X}(M)$ and $e \in F^{2k+1}$.
\end{lem}

\begin{proof}
Since $\bH^{2k}(M)/\bH^{2k}(M)_\Z \subset \G^0(M)$ acts on $\mathcal{V}(\lambda)$ via $\lambda$, the subspaces $\mathcal{V}(\lambda)$ and $\mathcal{V}(\lambda')$ are orthogonal for $\lambda \neq \lambda'$. So, if the action of $\eta \in \bH^{2k}(M)/\bH^{2k}(M)_\Z$ on an element $\phi \in \widehat{\oplus}_\lambda \mathcal{V}(\lambda)$ is the scalar multiplication of $e^{2\pi\im\lambda'(\eta)}$, then $\phi \in \mathcal{V}(\lambda')$. Now, for $\phi \in \mathcal{V}(\lambda)$, the action of $\eta$ on $\rho^\theta(\sigma(e)) \phi$ is the scalar multiplication of $e^{2\pi\im (\lambda + 2s(e))(\eta)}$, because $\theta$ commutes with the action of $\eta$. Hence $\rho^\theta(\sigma(e))\phi \in \mathcal{V}(\lambda + 2s(e))$, which leads to the lemma.
\end{proof}

Now, let $\lambda_i \in \mathcal{X}(M)$, ($i = 1, \ldots, 2^b$) be such that: $i \neq j$ implies $[ \lambda_i ] \neq [ \lambda_j ]$ in $\mathrm{Coker}\{ 2s : F^{2k+1} \to \mathcal{X}(M) \}$. Renumbering $\lambda_i$ such that $m(\lambda_i) > 0$, we obtain a subset $\{ \lambda_1, \ldots, \lambda_N \} \subset \mathcal{X}(M)$ with $N \le 2^b$. Then we can write:
$$
\widehat{\bigoplus}_{\lambda \in \mathcal{X}(M)} 
\mathcal{V}(\lambda) 
=
\bigoplus_{i = 1}^N \widehat{\bigoplus}_{\xi \in F^{2k+1}} 
\mathcal{V}(\lambda_i + 2s(\xi)).
$$
By Lemma \ref{lem:characterize_rep_free_part} and \ref{lem:Weyl_shift}, each subspace $\widehat{\bigoplus}_{\xi \in F^{2k+1}} \mathcal{V}(\lambda_i + 2s(\xi))$ is invariant under the action of $\G_\omega(M)$ through $\rho^\theta$. By Lemma \ref{lem:Weyl_shift}, we also have:
$$
\widehat{\bigoplus}_{\xi \in F^{2k+1}} 
\mathcal{V}(\lambda_i + 2s(\xi)) 
=
\widehat{\bigoplus}_{\xi \in F^{2k+1}} 
H_{\lambda_i + 2s(\xi)} \otimes \C^{m_i}
=
\H^\omega_{\lambda_i}|_{\G^0(M)} \otimes \C^{m_i},
$$
where we put $m_i = m(\lambda_i)$. Namely, $\G_\omega(M)$ also acts on $\widehat{\bigoplus}_{\xi \in F^{2k+1}} \mathcal{V}(\lambda_i + 2s(\xi))$ through $\rho^\omega_{\lambda_i} \otimes \id_{\C^{m_i}}$. As a result, there are two actions of $\G_\omega(M)$ on the same vector space. We know that these actions coincide on the subgroup $\G^0(M)$. The difference on $F^{2k+1}$ is expressed as follows:

\begin{prop} \label{prop:difference_shift}
There are unitary matrices $t^i_\xi \in U(m_i)$, $(\xi \in F^{2k+1})$ which make the following diagram commutative for all $e \in F^{2k+1}$:
$$
\begin{CD}
H_{\lambda_i + 2s(\xi)} \otimes \C^{m_i} 
@>{\id_{H_{\lambda_i + 2s(\xi)}} \otimes t^i_\xi}>>
H_{\lambda_i + 2s(\xi)} \otimes \C^{m_i} \\
@V{\rho^\omega_{\lambda_i}(\sigma(e)) \otimes \id_{\C^{m_i}}}VV
@VV{\rho^\theta(\sigma(e))}V \\
H_{\lambda_i + 2s(\xi + e)} \otimes \C^{m_i} 
@>>{\id_{H_{\lambda_i + 2s(\xi + e)}} \otimes t^i_{\xi + e}}>
H_{\lambda_i + 2s(\xi + e)} \otimes \C^{m_i}.
\end{CD}
$$
\end{prop}

\begin{proof}
By Lemma \ref{lem:Weyl_shift}, the restriction of $\rho^\theta(\sigma(e))$ provides the unitary isomorphism $H_{\lambda_i + 2s(\xi)} \otimes \C^{m_i} \to H_{\lambda_i + 2s(\xi + e)} \otimes \C^{m_i}$. By (\ref{formula:concrete_representation_action}), we also have the unitary isomorphism $\rho^\omega_{\lambda_i}(\sigma(e)) \otimes \id_{\C^{m_i}} : H_{\lambda_i + 2s(\xi)} \otimes \C^{m_i} \to H_{\lambda_i + 2s(\xi + e)} \otimes \C^{m_i}$. The unitary automorphism
$$
(\rho^\omega_{\lambda_i}(\sigma(e)) \otimes \id_{\C^{m_i}})^{-1} \circ
\rho^\theta(\sigma(e)) : \
H_{\lambda_i + 2s(\xi)} \otimes \C^{m_i} \longrightarrow 
H_{\lambda_i + 2s(\xi)} \otimes \C^{m_i}
$$
commutes with the action of $\G^0(M)$. Recall that $H_{\lambda_i}$ is irreducible. Thus, by Proposition \ref{prop:Schur_lemma}, there is a unitary matrix $u^i_\xi(e) \in U(m_i)$ which makes the following diagram commutative:
$$
\begin{CD}
H_{\lambda_i + 2s(\xi)} \otimes \C^{m_i} 
@=
H_{\lambda_i + 2s(\xi)} \otimes \C^{m_i} \\
@V{\rho^\theta(\sigma(e))}VV
@VV{\rho^\omega_{\lambda_i}(\sigma(e)) \otimes u^i_\xi(e)}V \\
H_{\lambda_i + 2s(\xi + e)} \otimes \C^{m_i} 
@=
H_{\lambda_i + 2s(\xi + e)} \otimes \C^{m_i}.
\end{CD}
$$
Since $\rho^\theta(\sigma(e))$ and $\rho^\omega_{\lambda_i}(\sigma(e))$ satisfy the relation in (\ref{formula:characterize_rep_free_part_II}), we obtain:
$$
u^i_\xi(e + e') = u^i_{\xi + e'}(e) u^i_\xi(e'). \quad
(\xi, e, e' \in F^{2k+1})
$$
Now, we define $t^i_\xi \in U(m_i)$ to be $t^i_\xi = u^i_0(\xi)$. Then the formula above gives $t^i_{e + \xi} = u^i_\xi(e) t^i_\xi$ for $\xi, e \in F^{2k+1}$, so the proposition is proved.
\end{proof}

\begin{cor} \label{cor:difference_shift}
There exists a unitary equivalence of projective representations of $\G_\omega(M)$:
$$
\Theta : \ 
\bigoplus_{i = 1}^N \left( \H^\omega_{\lambda_i} \otimes \C^{m_i} \right)
\longrightarrow
\H|_{\G_\omega(M)}.
$$
\end{cor}

\begin{proof}
The unitary matrices $t^i_\xi$ in Proposition \ref{prop:difference_shift} give the isometric map
$$
\oplus_\xi (\id_{H_{\lambda_i + 2s(\xi)}} \otimes t^i_\xi) : \
\bigoplus_{\xi \in F^{2k+1}} H_{\lambda_i + 2s(\xi)} \otimes \C^{m_i}
\longrightarrow
\bigoplus_{\xi \in F^{2k+1}} H_{\lambda_i + 2s(\xi)} \otimes \C^{m_i}.
$$
This extends to the unitary isomorphism:
$$
\mathcal{F}^i : \
\widehat{\bigoplus}_{\xi \in F^{2k+1}} \mathcal{V}(\lambda_i + 2s(\xi)) 
\longrightarrow
\widehat{\bigoplus}_{\xi \in F^{2k+1}} \mathcal{V}(\lambda_i + 2s(\xi))
$$ 
such that $\mathcal{F}^i \circ (\rho^\omega_{\lambda_i}(f) \otimes \id_{\C^{m_i}}) = \rho^\theta(f) \circ \mathcal{F}^i$ for $f \in \G_\omega(M)$. So the unitary equivalence $\Theta$ in this corollary is given by $\Theta = \theta \circ (\mathcal{F}^1 \oplus \cdots \oplus \mathcal{F}^N)$.
\end{proof}

\begin{thm} \label{thm:decomposition_admissible}
There are projective unitary representations $(\pi_i, V_i)$ of $T^{2k+1}$ with their cocycle $e^{2\pi\im L_M}$ and a unitary equivalence of projective representations of $\G(M)$:
$$
\Theta :  \ \bigoplus_{i = 1}^N \H^\omega_{\lambda_i, V_i}
\longrightarrow \H.
$$
\end{thm}

\begin{proof}
By Lemma \ref{lem:irreducible_free}, we have:
$$
\End_{\til{\G}_\omega(M)} 
(\oplus_{i = 1}^N \H^\omega_{\lambda_i} \otimes \C^{m_i}) =
\bigoplus_{i = 1}^N 
(\C \cdot \id_{\H^\omega_{\lambda_i}}) \otimes
\End(\C^{m_i}).
$$
The actions of $\G_\omega(M)$ and $T^{2k+1}$ on $\H$ commute with each other by means of Lemma \ref{lem:decompose_central_ext}. Since $\Theta$ in Corollary \ref{cor:difference_shift} is compatible with the action of $\G_\omega(M)$, there are, for $t \in T^{2k+1}$, unitary matrices $\pi_i(t) \in U(m_i)$ such that:
$$
\Theta^{-1} \circ \rho(t) \circ \Theta = 
\oplus_{i = 1}^N (\id_{\H^\omega_{\lambda_i}} \otimes \pi_i(t)).
$$
If we put $V_i = \C^{m_i}$, then $(\pi_i, V_i)$ is a projective unitary representation of $T^{2k+1}$ with its cocycle $e^{2\pi\im L_M}$. It is now direct to see that $\Theta$ gives rise to the unitary equivalence stated in the present theorem.
\end{proof}

We complete the proof of Theorem \ref{ithm:classification}: a finite dimensional projective representation of $T^{2k+1}$ is completely reducible \cite{Kar}. The irreducible decomposition of $(\pi_i, V_i)$ induces that of $\H^\omega_{\lambda_i, V_i}$ by Theorem \ref{thm:irreducibility}. Hence Theorem \ref{thm:decomposition_admissible} leads to Theorem \ref{ithm:classification} (a). Now Theorem \ref{ithm:classification} (b) follows from Corollary \ref{cor:irreducibility}.

\medskip

\begin{rem}
For a positive integer $\l$, we can also prove  Lemma \ref{lem:representation_identity_component} with $e^{2\pi\im S_M}$ replaced by $e^{2\pi\im \l S_M}$ in the same manner. So we can introduce the notion of admissible representations of $\G(M)$ of \textit{level} $\l$ by the same replacement in Definition \ref{idfn:admissible}. Almost the same argument leads to a classification similar to Theorem \ref{ithm:classification}: The difference is that the number of the equivalence classes of irreducible admissible representations of level $\l$ is $(2\l)^b r$, where $b = b_{2k+1}(M)$ is the Betti number, and $r$ is the number of the elements of the finite set 
$$
\{ t \in H^{2k+1}(M, \Z) |\ 2\l \cdot t = 0 \}.
$$
Note that, as in the case of $\l = 1$, admissible representations of $\G(S^1)$ of level $\l$ give rise to positive energy representations of $L\T$ of level $2\l$.
\end{rem}

\section{Properties of admissible representations}
\label{sec:properties}

This section is devoted to the proof of Theorem \ref{ithm:property}. Because of Theorem \ref{thm:decomposition_admissible}, it is enough to verify the properties for $(\rho^\omega_{\lambda, V}, \H^\omega_{\lambda, V})$. So $\omega$, $\lambda$ and $(\pi, V)$ are fixed in this section.

\subsection{Continuity}

We prove that the map $\rho^\omega_{\lambda, V} : \G(M) \times \H^\omega_{\lambda, V} \to \H^\omega_{\lambda, V}$ is continuous. For this aim, we remark the following fact, which can be seen readily: let $G$ be a topological group, and $(\rho, \H)$ a (projective) unitary representation of $G$ on a Hilbert space. Then $\rho : G \times \H \to \H$ is continuous, if and only if $\rho(\cdot)\phi : G \to \H$ is continuous for each $\phi \in \H$.

\begin{lem} \label{lem:continuity_free_part}
The map $\rho^\omega_\lambda : \G_\omega(M) \times \H^\omega_\lambda \to \H^\omega_\lambda$ is continuous.
\end{lem}

\begin{proof}
First, we show that $\rho^\omega_\lambda|_{\G^0(M)}(\cdot)\phi : \G^0(M) \to \H^\omega_\lambda$ is continuous for each $\phi \in \H^\omega_\lambda$. Since $\H^\omega_\lambda$ is the completion of $\oplus_{\xi \in F^{2k+1}} H_{\lambda + 2s(\xi)}$, it is enough to consider the case of $\phi \in \oplus_{\xi \in F^{2k+1}} H_{\lambda + 2s(\xi)}$. Then we can express $\phi$ as a finite sum $\phi = \phi_1 + \cdots + \phi_N$, where $\phi_i \in H_{\lambda + 2s(\xi_i)}$ for some $\xi_1, \ldots, \xi_N \in F^{2k+1}$. Thus $\rho^\omega_\lambda|_{\G^0(M)}(\cdot) \phi$ factors through $\oplus_{i = 1}^N H_{\lambda + 2s(\xi_i)}$, and is continuous by Lemma \ref{lem:representation_identity_component}. For $F^{2k+1}$ has the discrete topology, $\rho^\omega_\lambda(\sigma(\cdot))\phi : F^{2k+1} \to \H^\omega_\lambda$ is clearly continuous, where $\sigma$ is a harmonic splitting compatible with $\omega$. As is shown in Proposition \ref{prop:continuity_2_cocycle}, the group 2-cocycle $S_M$ is also continuous. Hence Lemma \ref{lem:characterize_rep_free_part} implies that $\rho^\omega_\lambda(\cdot)\phi : \G_\omega(M) \to \H^\omega_\lambda$ is continuous. Now, the remark at the beginning of this subsection completes the proof.
\end{proof}

\begin{thm}
The map $\rho^\omega_{\lambda, V} : \G(M) \times \H^\omega_{\lambda, V} \to \H^\omega_{\lambda, V}$ is continuous.
\end{thm}

\begin{proof}
It suffices to prove the continuity of $\rho^\omega_{\lambda, V}(\cdot)\phi : \G(M) \to \H^\omega_{\lambda, V}$ for each $\phi \in \H^\omega_{\lambda, V}$. This amounts to proving that the maps $\rho^\omega_{\lambda, V}(\cdot)\phi : \G_\omega(M) \to \H^\omega_{\lambda, V}$ and $\rho^\omega_{\lambda, V}(\sigma(\cdot))\phi : T^{2k+1} \to \H^\omega_{\lambda, V}$ are continuous for a harmonic splitting $\sigma$ compatible with $\omega$, because $S_M$ is continuous. The former is continuous by Lemma \ref{lem:continuity_free_part}, and so is the latter, because $T^{2k+1}$ has the discrete topology
\end{proof}

\subsection{Complexification}

As in \cite{Go}, we denote by $\G(M)_\C$ the $(2k+1)$th hypercohomology group of the complex of sheaves:
$$
\Z \longrightarrow
\u{A}^0_{M, \C} \overset{d}{\longrightarrow}
\u{A}^1_{M, \C} \overset{d}{\longrightarrow}
\cdots \overset{d}{\longrightarrow}
\u{A}^{2k}_{M, \C} \longrightarrow
0 \longrightarrow \cdots,
$$
where $\u{A}^q_{M, \C}$ is the sheaf of germs of $\C$-valued $q$-forms on $M$. We can think of $\G(M)_\C$ as a ``complexification'' of $\G(M)$. For example, if $M = S^1$, then $\G(S^1)_\C \cong L\C^*$, which is a complexification of $L\T$.

The properties of $\G(M)$ summarized in Section \ref{sec:Deligne_coh} hold for $\G(M)_\C$ under appropriate modifications, so that we have the natural exact sequence:
\begin{equation}
0 \to A^{2k}(M, \C)/A^{2k}(M, \C)_\Z \to \G(M)_\C \to H^{2k}(M, \Z) \to 0,
\label{formula:second_exact_complexification}
\end{equation}
where $A^{2k}(M, \C)$ is the group of $\C$-valued $2k$-forms on $M$, and $A^{2k}(M, \C)_\Z$ is the subgroup consisting of closed integral ones. The cup product and the integration yield the group 2-cocycle $S_M : \G(M)_\C \times \G(M)_\C \to \C/\Z$. For $\alpha, \beta \in \G^0(M)_\C = A^{2k}(M, \C)/A^{2k}(M, \C)_\Z$, we also have $S_{M}(\alpha, \beta) = \int_M \alpha \wedge d\beta \mod \Z$.

Since $A^{2k}(M, \C) = A^{2k}(M) \otimes \C$, the results in Section \ref{sec:topology} also hold for $\G(M)_\C$ under appropriate modifications: we have the isomorphism of topological groups
$$
\G^0(M)_\C \cong 
(\bH^{2k}(M, \C)/\bH^{2k}(M, \C)_\Z) \times d^*(A^{2k+1}(M, \C)),
$$
where $\bH^{2k}(M, \C)$ denotes the group of $\C$-valued harmonic $2k$-forms on $M$, and $\bH^{2k}(M, \C)_\Z = \bH^{2k}(M, \C) \cap A^{2k}(M, \C)_\Z$. Composing a harmonic splitting of (\ref{formula:second_exact_seq_again}) with the inclusion $\G(M) \to \G(M)_\C$, we get a harmonic splitting of (\ref{formula:second_exact_complexification}). Hence we can make $\G(M)_\C$ into a topological group such that $S_M$ is continuous. 

The properties of harmonic splittings of (\ref{formula:second_exact_seq_again}) in Section \ref{sec:harmonic_splitting} are carried to that of (\ref{formula:second_exact_complexification}) straightly. In particular, if a decomposition $\omega$ is fixed, then there is a canonical isomorphism $\G(M)_\C \cong \G^\omega(M)_\C \times T^{2k+1}$.

\begin{lem} \label{lem:complexification_identity_component}
There is an invariant dense subspace $E_\lambda \subset H_\lambda$ such that:
\begin{list}{}{\parsep=-2pt\topsep=4pt}
\item[(i)]
$(\rho_\lambda, E_\lambda)$ extends to a projective representation of $\G^0(M)_\C$;

\item[(ii)]
the map $\rho_\lambda(\cdot)\phi : \G^0(M)_\C \to E_\lambda$ is continuous for $\phi \in E_\lambda$.
\end{list}
\end{lem}

\begin{proof}
Recall Subsection \ref{subsec:Heisenberg} and Lemma \ref{lem:Heisenberg}: we have the dense subspace $E \subset H = \widehat{S}(W)$ and the projective representation $\rho : V_\C \times E \to E$ whose restriction $\rho : V \times E \to E$ extends to the projective unitary representation $\rho : V \times H \to H$. Thus, if we put $E_\lambda = E$, then $E_\lambda$ is clearly invariant under $\G^0(M)$. Note that $V_\C$ is the completion of $d^*(A^{2k+1}(M, \C))$, and $\lambda \in \mathcal{X}(M)$ extends to a homomorphism $\bH^{2k}(M, \C)/\bH^{2k}(M, \C)_\Z \to \C/\Z$. Hence we can extend $\rho_\lambda : \G^0(M) \times E_\lambda \to E_\lambda$ to a projective representation $\rho_\lambda : \G^0(M)_\C \times E_\lambda \to E_\lambda$. Since $\rho(\cdot)\phi : V_\C \to E$ is continuous, so is $\rho(\cdot)\phi : \G^0(M)_\C \to E_\lambda$.
\end{proof}

\begin{lem} \label{lem:complexification_free_part}
There is an invariant dense subspace $\E^\omega_\lambda \subset \H^\omega_\lambda$ such that:
\begin{list}{}{\parsep=-2pt\topsep=4pt}
\item[(i)]
$(\rho^\omega_\lambda, \E^\omega_\lambda)$ extends to a projective representation of $\G_\omega(M)_\C$;

\item[(ii)]
the map $\rho^\omega_\lambda(\cdot)\phi : \G_\omega(M)_\C \to \E^\omega_\lambda$ is continuous for $\phi \in \E^\omega_\lambda$.
\end{list}
\end{lem}

\begin{proof}
By means of the expression (\ref{formula:concrete_representation_space}), the subspace $\E^\omega_\lambda = \oplus_{\xi \in F^{2k+1}} \E_{\lambda + 2s(\xi)}$ is dense in $\H^\omega_\lambda$. We can see that the subspace is invariant under the action of $\G_\omega(M)$ by the help of the formula (\ref{formula:concrete_representation_action}). The formula and Lemma \ref{lem:complexification_identity_component} lead to (i). We can show (ii) in the same way as Lemma \ref{lem:continuity_free_part}.
\end{proof}

\begin{thm}
There is an invariant dense subspace $\E^\omega_{\lambda, V} \subset \H^\omega_{\lambda, V}$ such that:
\begin{list}{}{\parsep=-2pt\topsep=4pt}
\item[(i)]
$(\rho^\omega_{\lambda, V}, \E^\omega_{\lambda, V})$ extends to a projective representation of $\G(M)_\C$;

\item[(ii)]
the map $\rho^\omega_{\lambda, V}(\cdot)\phi : \G(M)_\C \to \E^\omega_{\lambda, V}$ is continuous for $\phi \in \E^\omega_{\lambda, V}$.
\end{list}
\end{thm}

\begin{proof}
The invariant dense subspace is given by $\E^\omega_{\lambda, V} = \E^\omega_\lambda \otimes V$. Then (i) and (ii) follow from Lemma \ref{lem:complexification_free_part}.
\end{proof}


\appendix

\section{Schur's lemma}

In this appendix, we show some general properties of representations on Hilbert spaces. First, we prove Proposition \ref{prop:Schur_lemma} (``Schur's lemma''):

\begin{prop} \label{appendix_prop:Schur_lemma}
Let $(\rho, \H)$ be a unitary representation of a group $G$ on a Hilbert space $\H$.

(a) If $\End_G(\H) = \C \cdot \id_{\H}$, then $(\rho, \H)$ is irreducible.

(b) Conversely, if $(\rho, \H)$ is irreducible, then $\End_G(\H) = \C \cdot \id_{\H}$.
\end{prop}

As in the main text, we mean by \textit{irreducible} that $\H$ contains no non-trivial invariant closed subspace.

\begin{proof}[The proof of Proposition \ref{appendix_prop:Schur_lemma} (a)]
Assume that $\H$ contains a non-trivial invariant closed subspace $\E$. By the \textit{projection theorem} \cite{R-S,Y}, we obtain the orthogonal decomposition $\H = \E \oplus \E^{\perp}$, where $\E^{\perp}$ is the orthogonal complement:
$$
\E^{\perp} = 
\{ v \in \H |\ (v, w) = 0 \ \mathrm{for} \ w \in \E \}.
$$
Let $P, P^{\perp} \in \End(\H)$ be the orthogonal projections onto $\E$ and $\E^{\perp}$, respectively. Since $\E^\perp$ as well as $\E$ is an invariant subspace, both $P$ and $P^\perp$ commute with the action of $G$, that is, $P, P^\perp \in \End_G(\H)$. By the assumption, $\E$ and $\E^\perp$ are non-trivial subspaces. Hence we have $P, P^\perp \not\in \C \cdot \id_{\H}$, so that $\End_G(\H) \supsetneq \C \cdot \id_{\H}$.
\end{proof}

To prove Proposition \ref{appendix_prop:Schur_lemma} (b), we make use of the \textit{spectral decomposition theorem} \cite{R-S,Y}. By definition, a family $\{ E(\lambda) \}_{\lambda \in \R}$ of orthogonal projections on the Hilbert space $\H$ is said to be a \textit{resolution of the identity} if it has the following properties:
\begin{align}
& \quad E(\lambda) E(\mu) = E(\mathrm{min}(\lambda, \mu)), 
\label{formula:monotonicity} \\
& \ \ \lim_{\varepsilon \downarrow 0}
\p{E(\lambda + \varepsilon)v - E(\lambda)v} = 0, \quad
(v \in \H) \\
& \lim_{\lambda \to - \infty} 
\p{E(\lambda)v} = 0, \
\lim_{\lambda \to + \infty}
\p{E(\lambda)v - v} = 0. \quad
(v \in \H)
\label{formula:limits}
\end{align}

\begin{prop}[The spectral decomposition theorem] 
There exists a one to one correspondence between self-adjoint bounded operators $A : \H \to \H$ and resolutions of the identity $\{ E(\lambda) \}_{\lambda \in \R}$ by the relation
\begin{equation} \label{formula:spectral_decomposition}
(Av, v) = 
\int \lambda d(E(\lambda)v, v). \quad (v \in \H)
\end{equation}
\end{prop}

For example, the multiplication of $c \in \R$ defines the self-adjoint operator $A = c \cdot \id_\H$. The resolution of the identity $\{ E(\lambda) \}$ corresponding to $A$ is
\begin{equation} \label{formula:scalar_multiplication_operator}
E(\lambda) =
\left\{
\begin{array}{cc}
\id_\H, & (\lambda \ge c) \\
0. & (\lambda < c)
\end{array}
\right.
\end{equation}

\begin{proof}[The proof of Proposition \ref{appendix_prop:Schur_lemma} (b)]
Notice that any $F \in \End(\H)$ admits the decomposition $F = \mathrm{Re} F + i \mathrm{Im} F$, where $\mathrm{Re} F$ and $\mathrm{Im} F$ are the self-adjoint operators given by
$$
\mathrm{Re} F = \frac{1}{2}(F + F^*), \quad
\mathrm{Im} F = \frac{1}{2i}(F - F^*).
$$
If $F$ commutes with the action of $G$, then so do $\mathrm{Re} F$ and $\mathrm{Im} F$. Thus, under the assumption in (b), it suffices to prove that: if $A \in \End_G(\H)$ is self-adjoint, then $A = c \cdot \id_\H$ for some $c \in \R$.

Let $\{ E_A(\lambda) \}$ denote the resolution of the identity corresponding to $A$ above. Then each $E_A(\lambda)$ also commutes with the action of $G$: to see this, we put $E'(\lambda) = \rho(g)^{-1} E_A(\lambda) \rho(g)$ for $g \in G$ fixed. We write $A'$ for the self-adjoint operator corresponding to the resolution of the identity $\{ E'(\lambda) \}$. Since $\rho(g)$ is unitary, we have $(E'(\lambda)v, v) = (E_A(\lambda)\rho(g)v, \rho(g)v)$. Thus (\ref{formula:spectral_decomposition}) yields
\begin{equation*}
\begin{split}
(A'v, v) 
&=
\int \lambda d(E'(\lambda) v, v) \\
&= 
\int \lambda d(E_A(\lambda)\rho(g)v, \rho(g)v)
=
(A\rho(g)v, \rho(g)v) = (Av, v).
\end{split}
\end{equation*}
Hence $A = A'$ and $E_A(\lambda) = E'(\lambda) = \rho(g)^{-1} E_A(\lambda) \rho(g)$.

Because of this commutativity, the closed subspaces in $\H$ given as the kernel and the range of $E_A(\lambda)$ are invariant under the action of $G$. Now, by the irreducibility of $\H$, these subspaces are $\{ 0 \}$ or $\H$, so that $E_A(\lambda)$ is $0$ or $\id_\H$. Then, by virtue of the properties (\ref{formula:monotonicity}) -- (\ref{formula:limits}), there uniquely exists $c \in \R$ that allows us to express $E_A(\lambda)$ as in (\ref{formula:scalar_multiplication_operator}). Therefore $A = c \cdot \id_\H$ as claimed.
\end{proof}

The next proposition indicates that there is no difference between \textit{unitary equivalence} and \textit{equivalence} in classifying unitary representations on Hilbert spaces.

\begin{prop} \label{appendix_prop:unitary_equivalence}
Let $(\rho_1, \H_1)$ and $(\rho_2, \H_2)$ be unitary representations of a group $G$ on Hilbert spaces. If there exists an equivalence between $\H_1$ and $\H_2$, then there also exists a unitary equivalence between them.
\end{prop}

\begin{proof}
We write $( \ , \ )_1$ and $( \ , \ )_2$ for the positive definite inner products on $\H_1$ and $\H_2$, respectively. Suppose that $\theta : \H_1 \to \H_2$ is an equivalence of the unitary representations. Let $( \ , \ )'_1$ be the positive definite inner product on $\H_1$ given by $(v, w)'_1 = (\theta v, \theta w)_2$. Since $\theta$ is an isomorphism, the inner products $( \ , \ )_1$ and $( \ , \ )'_1$ induce equivalent norms on $\H_1$. For $w \in \H_1$, we have the bounded operator $( \cdot, w )_1 : \H_1 \to \H_1$. Thus, applying the \textit{Riesz representation theorem} \cite{R-S,Y} to the Hilbert space $(\H_1, ( \ , \ )'_1)$, we obtain the bounded operator $A : \H_1 \to \H_1$ such that $(v, A w)'_1 = (v, w)_1$. We can readily see that $A$ is self-adjoint with respect to $( \ , \ )'_1$. Moreover, $A$ is invertible, positive and compatible with the action of $G$. Notice that $0$ does not belong to the spectrum $\sigma(A)$ by the positivity. Thus, as the functional calculus, we obtain the invertible self-adjoint bounded operator $\sqrt{A} : \H_1 \to \H_1$. In terms of the resolution of the identity $\{ E_A(\lambda) \}$ corresponding to $A$, we can express $\sqrt{A}$ as
$$
(\sqrt{A}v, v)'_1 = \int \sqrt{\lambda} d(E_A(\lambda)v, v)'_1.
$$
Because $A$ is compatible with the action of $G$, so is each $E_A$. Hence $\sqrt{A}$ is also compatible with the action of $G$. For $v, w \in \H_1$ we have
$$
(\theta \sqrt{A} v, \theta \sqrt{A} w)_2 = (\sqrt{A} v, \sqrt{A} w)'_1 = (v, A w)'_1 = (v, w)_1.
$$
Thus $\theta \circ \sqrt{A} : \ \H_1 \to \H_2$ is a unitary equivalence.
\end{proof}


\bigskip

\begin{acknowledgment}
Some of the ideas in Appendix were suggested by M. Furuta. I would like to express my gratitude to him. I would like to thank T. Tate for helpful comments on an earlier draft.
\end{acknowledgment}


\begin{flushleft}
Graduate school of Mathematical Sciences, University of Tokyo, \\
Komaba 3-8-1, Meguro-Ku, Tokyo, 153-8914 Japan. \\
e-mail: kgomi@ms.u-tokyo.ac.jp
\end{flushleft}

\end{document}